\documentclass[a4paper,11pt,reqno]{article}
\usepackage{latexsym, amsmath, amsfonts, amssymb, amsthm, amscd, epsfig}
\usepackage{url}

\usepackage[font=small, labelfont={sc}]{caption}

\usepackage[a4paper,scale={0.73,0.75},marginratio={1:1},footskip=10mm,headsep=10mm]{geometry}

\usepackage{color}

\usepackage{hyperref}

\numberwithin{equation}{section}

\def \beq {\begin{eqnarray}}
\def \eeq {\end{eqnarray}}
\def \beqn {\begin{eqnarray*}}
\def \eeqn {\end{eqnarray*}}

\newtheorem{theorem}{Theorem}[section]
\newtheorem{itlemma}[theorem]{Lemma}
\newtheorem{itproposition}[theorem]{Proposition}
\newtheorem{itcorollary}[theorem]{Corollary}
\newtheorem{itremark}[theorem]{Remark}
\newtheorem{itdefinition}[theorem]{Definition}
\newtheorem{itexample}[theorem]{Example}
\newtheorem{itclaim}[theorem]{Claim}
\newtheorem{itfact}[theorem]{Fact}
\newtheorem{itassumption}[theorem]{Assumption}

\newenvironment{fact}{\begin{itfact}\rm}{\end{itfact}}
\newenvironment{claim}{\begin{itclaim}\rm}{\end{itclaim}}
\newenvironment{lemma}{\begin{itlemma}}{\end{itlemma}}
\newenvironment{remark}{\begin{itremark}\rm}{\end{itremark}}
\newenvironment{corollary}{\begin{itcorollary}}{\end{itcorollary}}
\newenvironment{proposition}{\begin{itproposition}}{\end{itproposition}}
\newenvironment{definition}{\begin{itdefinition}\rm}{\end{itdefinition}}
\newenvironment{example}{\begin{itexample}\rm}{\end{itexample}}
\newenvironment{assumption}{\begin{itassumption}}{\end{itassumption}}

\newcommand{\be}[1]{\begin{equation}\label{#1}}
\newcommand{\ee}{\end{equation}}
\newcommand{\bl}[1]{\begin{lemma}\label{#1}}
\newcommand{\br}[1]{\begin{remark}\label{#1}}
\newcommand{\brs}[1]{\begin{remarks}\label{#1}}
\newcommand{\bt}[1]{\begin{theorem}\label{#1}}
\newcommand{\bd}[1]{\begin{definition}\label{#1}}
\newcommand{\bp}[1]{\begin{proposition}\label{#1}}
\newcommand{\bc}[1]{\begin{corollary}\label{#1}}
\newcommand{\bfact}[1]{\begin{fact}\label{#1}.}
\newcommand{\bex}[1]{\begin{example}\label{#1}.}
\newcommand{\ec}{\end{corollary}}
\newcommand{\efact}{\end{fact}}
\newcommand{\eex}{\end{example}}
\newcommand{\el}{\end{lemma}}
\newcommand{\er}{\end{remark}}
\newcommand{\ers}{\end{remarks}}
\newcommand{\et}{\end{theorem}}
\newcommand{\ed}{\end{definition}}
\newcommand{\ep}{\end{proposition}}
\newcommand{\epr}{\end{proof}}
\newcommand{\bpr}{\begin{proof}}
\newcommand{\bcl}[1]{\begin{claim}\label{#1}}
\newcommand{\ecl}{\end{claim}}
\newcommand{\bas}[1]{\begin{assumption}\label{#1}}
\newcommand{\eas}{\end{assumption}}

\newcommand{\ecs}{\end{corollary}}
\newcommand{\eers}{\end{exercise}}
\newcommand{\eexs}{\end{example}}
\newcommand{\eems}{\end{example}}
\newcommand{\els}{\end{lemma}}
\newcommand{\eles}{\end{lemmaex}}
\newcommand{\ets}{\end{theorem}}
\newcommand{\eds}{\end{definition}}
\newcommand{\eps}{\end{proposition}}
\newcommand{\bi}{\begin{itemize}}
\newcommand{\ei}{\end{itemize}}
\newcommand{\ben}{\begin{enumerate}}
\newcommand{\een}{\end{enumerate}}

\def\vbar{\mathchoice{\vrule height6.3ptdepth-.5ptwidth.8pt\kern-.8pt}
   {\vrule height6.3ptdepth-.5ptwidth.8pt\kern-.8pt}
   {\vrule height4.1ptdepth-.35ptwidth.6pt\kern-.6pt}
   {\vrule height3.1ptdepth-.25ptwidth.5pt\kern-.5pt}}
\def\fudge{\mathchoice{}{}{\mkern.5mu}{\mkern.8mu}}
\def\bbc#1#2{{\rm \mkern#2mu\vbar\mkern-#2mu#1}}
\def\bbb#1{{\rm I\mkern-3.5mu #1}}
\def\bba#1#2{{\rm #1\mkern-#2mu\fudge #1}}
\def\bb#1{{\count4=`#1 \advance\count4by-64 \ifcase\count4\or\bba A{11.5}\or
   \bbb B\or\bbc C{5}\or\bbb D\or\bbb E\or\bbb F \or\bbc G{5}\or\bbb H\or
   \bbb I\or\bbc J{3}\or\bbb K\or\bbb L \or\bbb M\or\bbb N\or\bbc O{5} \or
   \bbb P\or\bbc Q{5}\or\bbb R\or\bbc S{4.2}\or\bba T{10.5}\or\bbc U{5}\or
   \bba V{12}\or\bba W{16.5}\or\bba X{11}\or\bba Y{11.7}\or\bba Z{7.5}\fi}}

\makeatletter
\newtheorem*{rep@theorem}{\rep@title} \newcommand{\newreptheorem}[2]{%
\newenvironment{rep#1}[1]{%
\def\rep@title{\bf #2 \ref{##1} }%
\begin{rep@theorem} }%
{\end{rep@theorem} } }
\makeatother
\newreptheorem{theorem}{Theorem}
\newreptheorem{lemma}{Lemma}

\def\dim{\mathrm{dim}\,}
\def\Ker{\mathrm{Ker}\,}
\def\Id{\mathrm{Id}}
\def\Leb{\mathrm{Leb}\,}

\def \R {{\mathbb R}}

\def \PR {{\mathbb P}}
\def \E {{\mathbb E}}

\def \s {y}

\def \D{{\cal{D}}}

\newcommand{\ba}[1]{\addtocounter{for}{1} \begin{eqnarray}\label{#1}}
\newcommand{\ea}{\end{eqnarray}}

\def\sqr#1#2{{\vcenter{\vbox{\hrule height .#2pt
                             \hbox{\vrule width .#2pt height#1pt \kern#1pt
                                   \vrule width .#2pt}
                             \hrule height .#2pt}}}}

\def\pmb#1{\setbox0=\hbox{#1}%
   \kern-.025em\copy0\kern-\wd0
   \kern.05em\copy0\kern-\wd0
   \kern-.025em\raise.0433em\box0 }
\def\sqr#1#2{{\vcenter{\vbox{\hrule height.#2pt
     \hbox{\vrule width.#2pt height#1pt \kern#1pt
   \vrule width.#2pt}\hrule height.#2pt}}}}

\def\ve{\varepsilon}
\def\s{\sigma}
\def\d{\delta}
\def\l{\lambda}

\def\g{\gamma}

\def\a{\alpha}
\def\th{\theta}

\def\cal{\mathcal}

\newenvironment{myenumerate}{%
\begin{list}{\labelenumi}
	{%
	\setlength{\itemsep}{0.4em}%
	\setlength{\topsep}{0.5em}%
	\setlength\leftmargin{2.6em}%
	\setlength\labelwidth{2.15em}%
	\setlength{\labelsep}{0.45em}%
	\usecounter{enumi}%
	}%
	}%
{\end{list}
}

\renewenvironment{enumerate}{
\renewcommand{\theenumi}{\arabic{enumi}}%
\renewcommand{\labelenumi}{{\rm(\theenumi)}}%
\begin{myenumerate}}%
{\end{myenumerate}}

{\end{myenumerate}}

{\end{myenumerate}}

\newenvironment{myitemize}{%
\begin{list}{$\bullet$}%
 	{%
	\setlength{\itemsep}{0.4em}%
	\setlength{\topsep}{0.5em}%
	\setlength\leftmargin{2.6em}%
	\setlength\labelwidth{2.15em}%
	\setlength{\labelsep}{0.45em}%
	}%
	}%
{\end{list}}

\renewenvironment{itemize}{
\begin{myitemize}}%
{\end{myitemize}}

\begin{document}

\parindent 0pt

\title{Diffusions under a local strong H\"ormander condition.\\Part II: tube estimates}
\author{ \textsc{Vlad Bally}\thanks{%
Universit\'e Paris-Est, LAMA (UMR CNRS, UPEMLV, UPEC), MathRisk INRIA, F-77454
Marne-la-Vall\'{e}e, France. Email: \texttt{bally@univ-mlv.fr} }\smallskip \\
%EndAName
\textsc{Lucia Caramellino}\thanks{%
Dipartimento di Matematica, Universit\`a di Roma - Tor Vergata, Via della
Ricerca Scientifica 1, I-00133 Roma, Italy. Email: \texttt{%
caramell@mat.uniroma2.it}. }\smallskip\\
\textsc{Paolo Pigato}\thanks{%
INRIA, Villers-l\`es-Nancy, F-54600, France
Universit\'e de Lorraine, IECL, UMR 7502, Vandoeuvre-l\`es-Nancy, F-54600, France \texttt{%
paolo.pigato@inria.fr}. }\smallskip\\
}
\maketitle

\begin{abstract}
We study lower and upper bounds for the probability that a diffusion process in $\R^n$ remains in a tube around a skeleton path up to a fixed time. We assume that the diffusion coefficients $\sigma_1,\ldots,\sigma_d$ may degenerate but they satisfy a strong H\"ormander condition involving the first order Lie brackets around the skeleton of interest. The tube is written in terms of a norm which accounts for the non-isotropic structure of the problem: in a small time $\delta$, the diffusion process propagates with speed $\sqrt{\delta}$ in the direction of the diffusion vector fields $\sigma _{j}$ and with speed $\delta=\sqrt{\delta}\times \sqrt{\delta}$ in the direction of $[\sigma _{i},\sigma _{j}]$. The proof consists in a concatenation technique which strongly uses the lower and upper bounds for the density proved in the part I.
\end{abstract}

\tableofcontents

\date{}

\section{Introduction}
We consider a diffusion process in $\R^n$ solution of
\[
dX_{t}=\sum_{j=1}^{d}\s
_{j}(t,X_{t})\circ dW_{t}^{j}+b(t,X_{t})dt,\quad \quad X_0=x_0.
\]
where $W=(W^{1},...,W^{d})$ is a standard Brownian motion and $\circ
dW_{t}^{j} $ denotes the Stratonovich integral.
We assume suitable regularity properties for $\s_{j},b:\R^{+}\times \R^{n}\rightarrow \R^{n}$ (see \eqref{Not2} for details). We also assume that the coefficients $\sigma_{j},b$ verify the strong H\"{o}rmander condition of order one (that is, involving the $\sigma _{j}$'s and their first order Lie brackets $[\sigma _{i},\sigma_{j}]$'s) locally around a skeleton path
\[
dx_{t}(\phi )=\sum_{j=1}^{d}\sigma _{j}(t,x_{t}(\phi ))\phi
_{t}^{j}dt+b(t,x_{t}(\phi ))dt,\quad \quad x_0(\phi)=x_0
\]
(this is formally written in property $(H_2)$ of \eqref{Not6}).
In such a framework, in this paper we find exponential lower and upper bounds for the probability that the diffusion $X$ remains in a small tube around the skeleton path $x(\phi)$.

Several works have considered this subject, starting from Stroock and Varadhan in \cite{StroockVaradhan:72}, where such result is used to prove the support theorem for diffusion processes. In their work, the tube is written in terms of the Euclidean norm, but later on different norms have been used to take into account the regularity of the trajectories (\cite{BenArousGradinaruLedoux:1994, FrizLyonsStrook:06}) and their geometric structure (\cite{pigato:14}). This kind of problems is also related to the Onsager-Machlup functional and large or moderate deviation theory, see e.g. \cite{capitaine:00,IkedaWatanabe:89,guillin2003}.

In this work, we construct the tube using a distance coming from a norm which reflects the non isotropic structure of the problem, i.e. the fact that the diffusion process $X_{t}$ propagates with speed $\sqrt{t}$ in the direction of the diffusion vector fields $\sigma _{j}$ and with speed $t=\sqrt{t}\times \sqrt{t}$ in the direction of $[\sigma _{i},\sigma _{j}]$. We also prove that this distance  is locally equivalent with the standard control (Carath\'eodory) metric.

A key step in proving our tube estimates is given by the use of the density estimates provided in \cite{BCP1}. Generally speaking, there is a strong connection between tube and density estimates. In this work we use a concatenation of short time density estimates to prove a tube estimate, but one may proceed in reverse order:
tubes estimates, for instance, can be used to provide lower bounds for the density. In \cite{[BFM]}, tube estimates for locally elliptic diffusions are proved, and applied to find lower bounds for the probability to be in a ball at fixed time and bounds for the distribution function. In \cite{BDM}, this is applied to lognormal-like stochastic volatility models, finding estimates for the tails of the distribution, and estimates on the implied volatility. 

The paper is organized as follows. In Section \ref{notationsresults}, we state our main result, given in Theorem \ref{mttubesS}, and we propose some examples of application. The proof of Theorem \ref{mttubesS} is developed in Section \ref{diffusion}. In Section \ref{sectioncontrol} we study the local equivalence between the control metric and the distance we use to  define the tube when the diffusion coefficients depend on the space variable only. As a straightforward consequence, we can state our tube estimate result in terms of the Carath\'eodory metric (see Theorem \ref{tubecontrold}).

\section{Notation and main results}\label{notationsresults}
We recall the notation from \cite{BCP1} and introduce some new ones. We consider vector fields $\s_{j},b:\R^{+}\times \R^{n}\rightarrow \R^{n}$ which are four time differentiable in $x\in \R^{n}$ and one time differentiable in time $t\in \R^{+}$, and suppose that the derivatives with respect to the space $x\in \R^{n}$ are one time differentiable with respect to $t$.

Hereafter, for $k\geq 1$, $\alpha =(\alpha _{1},\ldots,\alpha _{k})\in \{1,...,n\}^{k}$ represents a multi-index with length $|\alpha| =k$ and $%
\partial _{x}^{\alpha }=\partial _{x_{\alpha _{1}}}\cdots\partial _{x_{\alpha
_{k}}}$. We allow the case $k=0$ by setting $\alpha=\emptyset$ (the void multiindex), $|\alpha|=0$ and $\partial^\alpha_x=\Id$.

For $(t,x)\in \R^{+}\times \R^{n}$ we denote by $n(t,x)$ a constant such that
%for every $%
%s\in [ (t-1)\vee 0,t+1],y\in B(x,1)$ and for every multi index $\alpha $ of length less than or equal to four
%\be{Not2}
%\left\vert \partial _{x}^{\alpha }b(s,y)\right\vert +\left\vert \partial
%_{t}\partial _{x}^{\alpha }b(s,y)\right\vert +\sum_{j=1}^{d}\left\vert
%\partial _{x}^{\alpha }\sigma _{j}(s,y)\right\vert +\left\vert \partial
%_{t}\partial _{x}^{\alpha }\sigma _{j}(s,y)\right\vert )\leq n(t,x).
%\ee
%
%\be{Not2}
%\begin{array}{rl}
%\displaystyle
%\sup_{s\in[ (t-1)\vee 0,t+1]}\sup_{y\in B(x,1)}\sup_{|\alpha|\leq 4}
%&\Big(\left\vert \partial _{x}^{\alpha }b(s,y)\right\vert +\left\vert \partial
%_{t}\partial _{x}^{\alpha }b(s,y)\right\vert  \smallskip\\
%\displaystyle
%&\ +\sum_{j=1}^{d}\left\vert
%\partial _{x}^{\alpha }\sigma _{j}(s,y)\right\vert+\left\vert \partial
%_{t}\partial _{x}^{\alpha }\sigma _{j}(s,y)\right\vert \Big)
%\leq n(t,x).
%\end{array}
%\ee
%
\be{Not2}
\begin{array}{c}
\forall s\in[ (t-1)\vee 0,t+1],\ \forall y\in B(x,1)\mbox{ one has }\smallskip\\
\displaystyle
\sum_{|\alpha|=0}^4\Big(\left\vert \partial _{x}^{\alpha }b(s,y)\right\vert +\left\vert \partial
_{t}\partial _{x}^{\alpha }b(s,y)\right\vert  +\sum_{j=1}^{d}\left\vert
\partial _{x}^{\alpha }\sigma _{j}(s,y)\right\vert+\left\vert \partial
_{t}\partial _{x}^{\alpha }\sigma _{j}(s,y)\right\vert \Big)
\leq n(t,x).
\end{array}
\ee

For $f,g:\R^{+}\times \R^{n}\rightarrow \R^{n}$ we define the directional derivative (w.r.t. the space variable $x$) $\partial_{g}f(t,x)=\sum_{i=1}^{n}g^{i}(t,x)\partial _{x_{i}}f(t,x)$, and  we recall that the Lie bracket (again w.r.t. the space variable)
is defined as $[g,f](t,x)=\partial _{g}f(t,x)-\partial _{f}g(t,x)$.
Let $M\in \mathcal{M}_{n\times m}$ be a matrix with full row rank. We write $M^T$ for the transposed matrix, and $MM^T$ is invertible. We denote by $\l_*(M)$ (respectively $\l^*(M))$ the smallest (respectively the largest) singular value of $M$.
We recall that singular values are the square roots of the eigenvalues of $MM^T$, and that, when $M$ is symmetric, singular values coincide with the absolute values of the eigenvalues of $M$. In particular, when $M$ is a covariance matrix, $\l_*(M)$ and $\l^*(M)$ coincide with the smallest and the largest eigenvalues of $M$.

We consider the following norm on $\R^{n}$:
\be{Not4}
\left\vert y\right\vert _{M}=\sqrt{\left\langle (MM^T)^{-1}y,y\right\rangle}.
\ee
We introduce the $n\times d^2$ matrix $A(t,x)$ defined as follows. We set $m=d^2$ and define the function
\be{lip}
l(i,p)=(p-1)d+i\in\{1,\dots,m\},\quad p,i\in\{1,\dots,d\}.
\ee
Notice that $l(i,p)$ is invertible.  For $l=1,\ldots,m$, we set the (column) vector field $A_{l}(t,x)$ in $\R^n$ as follows:
\be{Al}
\begin{split}
A_{l}(t,x) &=[\s_{i},\s_{p}](t,x) \quad \mbox{if}\quad l=l(i,p) \quad \mbox{with}\quad i\neq p, \\
&= \s_{i}(t,x)\quad \mbox{if}\quad l=l(i,p) \quad \mbox{with}\quad i= p
\end{split}
\ee
and we set the $n\times m$ matrix $A(t,x)$ to be the one having $A_1(t,x),\ldots,A_m(t,x)$
as its columns, that is
\be{Alucia}
A(t,x)=[A_1(t,x),\ldots,A_m(t,x)].
\ee
We denote by $\lambda (t,x)$ the smallest singular value of $A(t,x)$, so
\be{Not9'}
\lambda (t,x)^2=
\l_*(A(t,x))^2= \inf_{\left\vert \xi \right\vert =1}\sum_{i=1}^{m}\left\langle
A_{i}(t,x),\xi \right\rangle ^{2}.
\ee
For fixed $R>0$ we define the $m\times m$ diagonal scaling matrix $D_R$ as
\begin{equation}\label{DR}
\begin{split}
(D_R)_{l,l}&=R \quad \mbox{if}\quad l=l(i,p) \quad \mbox{with}\quad i\neq p, \\
&= \sqrt{R} \quad \mbox{if}\quad l=l(i,p) \quad \mbox{with}\quad i= p
\end{split}
\end{equation}
and the scaled directional matrix
\be{Adlucia}
A_R(t,x)=A(t,x)D_R.
\ee
Notice that the $l$th column of the matrix $A_R(t,x)$ is given by
$\sqrt{R}\s_{i}(t,x)$ if $l=l(i,p)$ with $ i= p$, and if $i\neq p$ then the $l$th column of $A_R(t,x)$ is $R[\s_{i},\s_{p}](t,x)=[\sqrt{R}\s_{i},\sqrt{R}\s_{p}](t,x)$.

For a control $\phi \in L^2([0,T],\R^n)$ we consider the skeleton $x(\phi)$ associated to \eqref{equation}, that is,
\be{skeleton}
dx_{t}(\phi )=\sum_{j=1}^{d}\sigma _{j}(t,x_{t}(\phi ))\phi
_{t}^{j}dt+b(t,x_{t}(\phi ))dt, \quad x_0(\phi)=x_0.
\ee
In the following, we also need a function $R\,:\, [0,T]\rightarrow (0,1]$ that will play the role of a radius function (for the tube around $x(\phi)$).

We consider now a ``regularity property''  already introduced in \cite{BallyKohatsu:10}, which  is needed to control the growth of certain quantities along the skeleton path.
For $\mu \geq 1$ and $0<h\leq 1$ we denote by $L(\mu ,h)$ the following class
of functions:
\begin{equation}\label{Lhmu}
L(\mu,h)=\big\{f:\R^{+}\rightarrow \R^{+}\mbox{ such that }f(t)\leq \mu f(s)\quad \mbox{for}\quad \left\vert t-s\right\vert \leq h\big\}.
\end{equation}

From now on, we make use of the following hypotheses: there exist some functions $%
n:[0,T]\rightarrow \lbrack 1,\infty )$ and $\lambda :[0,T]\rightarrow (0,1]$
such that for some $\mu \geq 1$ and $0<h\leq 1$ we have
\begin{equation}
\begin{array}{lcl}
(H_{1})\quad & n(t,x_{t}(\phi ))  \leq n_{t},\quad\forall t\in \lbrack
0,T],\smallskip \\
(H_{2})\quad & \lambda (t,x_{t}(\phi ))  \geq \lambda _{t},\quad\forall t\in
\lbrack 0,T],\smallskip \\
(H_{3})\quad & R_.,\,\left\vert \phi _{.}\right\vert ^{2},\,n_{.},\,\lambda _{.}  \in
L(\mu ,h).
\end{array}
\label{Not6}
\end{equation}
Recall that $\phi \in L^2([0,T],\R^n)$ is the control giving the skeleton path and $R\,:[0,T]\to (0,1]$ stands for the radius function.

\br{dghd}
Hypothesis $(H_{2})$ implies that for each $t\in (0,T),$ the space $\R^n$
is spanned by the vectors $(\sigma _{i}(t,x_{t}),[\sigma _{j},\sigma
_{p}](t,x_{t}))_{i,j,p=1,...,d,j<p}$, meaning that a strong H\"{o}rmander condition locally holds along the curve $x_{t}(\phi ).$
\er
Let $X$ denote a process in $\R^n$  solving
\be{equation}
dX_{t}=\sum_{j=1}^{d}\s
_{j}(t,X_{t})\circ dW_{t}^{j}+b(t,X_{t})dt,\quad \quad X_0=x_0,
\ee
$W$ being  a standard Brownian motion in $\R^d$.
Remark that $(H_1)$ is only a local assumption: we do not assume global Lipschitz continuity or sublinear growth properties  for the coefficients, so the above SDE might not have a unique solution.
We only assume to work with a continuous adapted process $X$ solving \eqref{equation} on the time interval $[0,T]$.

For $K,q,K_*,q_*>0$, $\mu \geq 1$, $h\in(0, 1]$, $n\,:\,[0,T]\to [1,+\infty)$, $\lambda\,:\,[0,T]\to (0,1]$ and $\phi\in L^2([0,T],\R^n)$,  we set the functions
\begin{equation}\label{H-R}
\begin{array}{l}
H_t=K \left(\frac{\mu n_t}{\l_t}\right)^q,\smallskip\\
R_t^*(\phi)=\exp\left(-K_* \left(\frac{\mu n_t}{\l_t}\right)^{q_*} \mu^{2q_*}\right)
\left( h\wedge \inf_{0\leq \d\leq h}
\left\{
\d \big/ \int_t^{t+\d} |\phi_s|^2 ds \right\}
\right).
\end{array}
\end{equation}
The main result of this paper is the following:
\begin{theorem}\label{mttubesS}
Let $\mu\geq 1$, $h\in(0,1]$, $n\,:\,[0,T]\to [1,+\infty)$, $\lambda\,:\,[0,T]\to (0,1]$, $R\,:\,[0,T]\to (0,1]$ and $\phi\in L^2([0,T],\R^n)$ be such that $(H_1)$--$(H_3)$ in \eqref{Not6} hold.
Then there exist $K,q,K_*,q_*>0$ such that, for $H$ and $R^*(\phi)$ as in \eqref{H-R}, if $R_t\leq R_t^*(\phi)$ one has
\begin{equation}\label{tuberesultstrong}
\begin{array}{rl}
\displaystyle
\exp\left(- \int_0^T H_t \left(\frac{1}{R_t}+|\phi_t|^2\right)dt \right)
\leq
&\displaystyle
\PR\left(\sup_{t\leq T} |X_t-x_t(\phi)|_{A_{R_t}(t,x_t(\phi))}\leq 1 \right)\\
&\displaystyle
\leq \exp\left(- \int_0^T e^{-H_t} \left(\frac{1}{R_t}+|\phi_t|^2\right)dt \right).
\end{array}
\end{equation}
\end{theorem}
The proof of Theorem \ref{mttubesS} is developed in Section \ref{diffusion}. We discuss here some comments and examples.

\br{1}
The estimate \eqref{tuberesultstrong} allows for a regime shift, meaning that the dimension of the space generated by the $\s_i$'s and the $[\sigma_i,\sigma_j]$'s may change along the tube, and this is accounted by the variation of $A_R$ along $x_t(\phi)$.
\er
\br{4}
The fact that $R\in L(\mu,h)$ implies that $\inf_{t\in [0,T]} R_t>0$. So, the radius of the tube is small, but cannot go to $0$ at any time.
\er
\begin{remark}
The lower bound holds even if the inequality $R_t\leq R_t^*(\phi)$ is not satisfied, in the form
\[
\exp\left(- \int_0^T H_t \left(\frac{1}{h}+\frac{1}{R_t}+|\phi_t|^2\right)dt \right)
\leq\PR\left(\sup_{t\leq T} |X_t-x_t(\phi)|_{A_{R_t}(t,x_t(\phi))}\leq 1 \right).
\]
Details are given in next Theorem \ref{th-lower}.
\er
\br{jfgsjgt}
Suppose $X_{t}=W_{t}$ and $x(\phi)=0$, so that $n_{t}=1$, $\l_{t}=1$, $\mu =1$ and $\phi _{t}=0$. Take $R_t=R$ constant. Then $%
\left\vert X_{t}-x_{t}(\phi)\right\vert _{A_{R}(t,x_{t}(\phi ))}=R^{-1/2}W_{t}$ and
we obtain $\exp (-C_1 T/R)\leq \PR(\sup_{t\leq T}\left\vert W_{t}\right\vert \leq \sqrt{R})\leq \exp (-C_2 T/R)$ which is consistent with the standard estimate (see
\cite{IkedaWatanabe:89}).
\er

A global two-sided bound for the density of $X_t$ is proved in \cite{KusuokaStroock:87}, under the \emph{strong} H\"{o}rmander non-degeneracy condition. It is also assumed that the coefficients do not depend time, i.e. $b(t,x)=b(x),\, \s(t,x)=\s(x)$, and that $b(x)=\sum_{j=1}^d \a_i \s_i(x)$, with $\a_i\in C_b^\infty(\R^n)$ (i.e. the drift is generated by the vector fields of the diffusive part, which is a quite restrictive hypothesis). This bound is Gaussian in the control metric that we now define. For $x,y\in \R^{n}$ we denote by $C(x,y)$
the set of controls $\psi \in L^{2}([0,1];\R^d)$ such that the
corresponding solution of
\[
du_{t}(\psi )=\sum_{j=1}^{d}\sigma _{j}(u_{t}(\psi))\psi _{t}^{j}dt, \quad\quad u_{0}(\psi )=x
\]
satisfies $u_{1}(\psi )=y$. The control (Carath\'eodory) distance is defined as
\begin{equation*}
d_{c}(x,y)=\inf \Big\{\Big(\int_{0}^{1}\left\vert \psi _{s}\right\vert ^{2}ds%
\Big)^{1/2}:\psi \in C(x,y)\Big\}.
\end{equation*}
The result in \cite{KusuokaStroock:87} is the following. Let $p_\delta(x,\cdot)$ denote the density of $X_\delta$ with starting condition $X_0=x$. Then there exists a constant $M\geq 1$ such that
\[
\begin{split}
&\frac{1}{M |B_{d_c}(x,\sqrt \delta)|}\exp\left( -\frac{M d_c(x,y)^2}{\delta}\right)\\
&\quad \quad \quad \leq p_\delta(x,y)
\leq \frac{M}{|B_{d_c}(x,\sqrt\delta)|}\exp\left( -\frac{d_c(x,y)^2}{M \delta}\right)
\end{split}
\]
where $\delta\in(0,T]$, $x,y\in\R^n$, $B_d(x,r)=\{y\in \R^n: d(x,y)<r\}$ and $|B_{d_c}(x,r)|$ denotes its Lebesgue measure. Remark that now, as in \cite{KusuokaStroock:87}, $\s(t,x)=\s(x)$. We define the semi distance $d$ via: $d(x,y) < \sqrt{R}$ if $|x-y|_{A_R(x)}<1$, and prove in Section \ref{sectioncontrol} the local equivalence of $d$ and $d_c$. This allows us to state Theorem \ref{mttubesS} in the control metric:
\bt{tubecontrold}
Suppose that the diffusion  coefficients  $\sigma_j$, $j=1,\ldots,d$, in \eqref{equation} depend on the space variable $x$ only and that the hypotheses of Theorem \ref{mttubesS} hold. Then,
\begin{multline}
\exp\left(- \int_0^T H_t \left(\frac{1}{R_t}+|\phi_t|^2\right)dt \right)\\
\leq \PR\left(\sup_{0\leq t\leq T} d_c(X_t,x_t(\phi))\leq \sqrt{R_t} \right)\leq
\exp\left(- \int_0^T e^{-H_t} \left(\frac{1}{R_t}+|\phi_t|^2\right)dt \right).
\end{multline}
\et

We prove the tube estimates in Section \ref{diffusion}, whereas the equivalence between the matrix norm and the Carath\'eodory distance is given in Section \ref{sectioncontrol}.

\medskip

We present now two examples of application.

\medskip

\textbf{Example 1}. {\small \textsc{[Grushin diffusion]}} Consider a positive, fixed $R$ and the two dimensional diffusion process
\begin{equation*}
X_{t}^{1}=x_{1}+W_{t}^{1},\quad
X_{t}^{2}=x_{2}+\int_{0}^{t}X_{s}^{1}dW_{s}^{2}.
\end{equation*}%
%We have $\s_1(x)=\left(\begin{array}{cc} 1 \\ 0    \end{array} \right)$, $\s_2(x)=\left(\begin{array}{cc} 0 \\ x_1     \end{array} \right)$, $[\s_1,\s_2](x)=\left(\begin{array}{cc} 0 \\ 1 \end{array} \right)$.
Here
$$
A_R A_R^T (x)=\left(\begin{array}{cc}
R & 0 \\
0 & R(x_1^2+2R)
\end{array}  \right),
$$
so the associated norm is
$|\xi|_{A_R(x)}^2= \frac{\xi_1^2}{R}+\frac{\xi_2^2}{R(x_1^2+2R)}$.
On $\{x_1= 0\}$, $|\xi|_{A_R(x)}^2= \frac{\xi_1^2}{R}+\frac{\xi_2^2}{2R^2}$
and consequently $\{\xi :\left\vert \xi \right\vert _{A_{R }(x)}\leq 1\}$ is an ellipsoid.

If we take a path $x(t)$ with $x_1(t)$ which keeps far from zero then we have ellipticity along the path and we may use estimates for elliptic SDEs (see \cite{[BFM]}). If $x_{1}(t)=0$ for some $t\in \lbrack 0,T]$
we need our estimate. Let us compare the norm in the two cases: if $x_{1}>0$ the diffusion matrix is non-degenerate and
we can consider the norm $\left\vert \xi \right\vert _{B_{R}(x)}$ with
$B_R(x)=R \sigma(x).$ We have%
\begin{equation*}
\left\vert \xi \right\vert _{B_{R }(x)}^{2}=\frac{1}{R }\xi
_{1}^{2}+\frac{1}{R x_{1}^{2}}\xi _{2}^{2}\geq \frac{1}{R }\xi
_{1}^{2}+\frac{1}{R (x_1^2+2R)}\xi _{2}^{2}=\left\vert \xi
\right\vert _{A_{R }(x)}^{2},
\end{equation*}%
and the two norms are equivalent for $R$ small.
Let us now take $x_{t}(\phi )=(0,0).$ We have $n_{s}=1$ and
$\lambda _{s}=1$ and $X_{t}-x_{t}(\phi)=(W_{t}^{1},%
\int_{0}^{t}W_{s}^{1}dW_{s}^{2})$, so we obtain%
\begin{equation*}
\begin{split}
e^{-C_1 T/R }&\leq
\PR\left(\sup_{t\leq T}\left\{ \frac{1}{R}\left\vert W_{t}^{1}\right\vert ^{2}+%
\frac{1}{2R^{2}}\left\vert \int_{0}^{t}W_{s}^{1}dW_{s}^{2}\right\vert
^{2}\right\} \leq 1\right)\\
&=\PR\left(\sup_{t\leq T}(\left\vert X_{t}-x_{t}\right\vert
_{A_R(x_t)}^{2}\leq 1 \right)\leq e^{-C_2 T/R }.
\end{split}
\end{equation*}

\noindent
\textbf{Example 2.} {\small \textsc{[Principal invariant diffusion on the Heisenberg group]}}
Consider on $\R^3$ the vector fields $\partial_{x_1}-\frac{x_2}{2} \partial_{x_3}$ and $\partial_{x_2}-\frac{x_1}{2} \partial_{x_3}$. The associated Markov process is the triple given by a Brownian motion on $\R^2$ and its L\'evy area, that is
\[
X_{t}^{1}=x_{1}+W_{t}^{1},\quad X_{t}^{2}=x_{2}+W_{t}^{2},\quad
X_{t}^{3}=x_{3}+\frac{1}{2}\int_{0}^{t}X_{s}^{1}dW_{s}^{2}-
\frac{1}{2} \int_{0}^{t}X_{s}^{2}dW_{s}^{1}.
\]
We refer e.g. to \cite{Driver2005340,Bakry20081905,Li2007497}, where gradient bounds for the heat kernel are obtained, and \cite{bau}.
Since the diffusion is in dimension $n=3$ and the driving Brownian in dimension $d=2$, ellipticity cannot hold. Direct computations give
\[
\s_1(x)=
\left(\begin{array}{c}
1 \\
0 \\
- \frac{x_2}2
\end{array}\right), \quad
\s_2(x)=
\left(\begin{array}{c}
0 \\
1 \\
 \frac{x_1}2
\end{array}\right), \quad
[\s_1,\s_2](x)=
\partial_{\s_1}\s_2-
\partial_{\s_2}\s_1
=
\left(\begin{array}{c}
0 \\
0 \\
1
\end{array}\right).
\]
Therefore $\s_1(x),\s_2(x),[\s_1,\s_2](x)$ span $\R^3$ and  hypoellipticity holds.
In $x=0$ we have $|\xi|^2_{A_R(0)}= \frac{\xi_1^2+\xi_2^2}{R} + \frac{\xi_3^2}{2 R^2}$, so taking the control $\phi\equiv 0$ and denoting  $A_{t}(W)=\frac{1}{2}\int_{0}^{t}X_{s}^{1}dW_{s}^{2}-
\frac{1}{2} \int_{0}^{t}X_{s}^{2}dW_{s}^{1}$ (the L\'evy area), we obtain
\[
\begin{split}
\PR\left( \sup_{t\leq T/R} |W^1_t|^2+|W^2_t|^2 + \frac{|A_t(W)|^2}{2} \leq 1 \right)
&=
\PR\left( \sup_{t\leq T} \frac{|W^1_t|^2+|W^2_t|^2}{R} + \frac{|A_t(W)|^2}{2 R^2} \leq 1 \right)\\
&=
\PR\left(\sup_{t\leq T} |X_t|^2_{A_R(x_t(\phi))}\leq 1	\right).
\end{split}
\]
Appling our estimate we have
\[
e^{-C_1 T/R}\leq
\PR\left( \sup_{t\leq T/R} |W^1_t|^2+|W^2_t|^2 + \frac{|A_t(W)|^2}{2} \leq 1 \right) \leq e^{-C_2 T/R}.
\]

\section{Tube estimates}\label{diffusion}
The proof of Theorem \ref{mttubesS} is inspired by the  approach in \cite{[BFM]}. A similar procedure is also used in \cite{pigato:14} in a weak H\"ormander framework. Such a proof strongly uses the estimates for the density developed in \cite{BCP1} and it is crucial that these estimates hold in a time interval of a fixed small length. This is because the proof consists in a ``concatenation'' of such estimates in order to recover the whole time interval $[0,T]$. And since the ``concatenation'' works around the skeleton path $x(\phi)$, it suffices that the properties for all objects hold only locally around $x(\phi)$, as required in \eqref{Not6}. In order to set-up this program, we need the precise behavior of the norm $|\cdot |_{A_R}$. So, we first present the desired properties for $|\cdot |_{A_R}$ (Section \ref{norms}) and then we proceed with the proof of Theorem \ref{mttubesS} (Section \ref{diffusion}).

\subsection{Matrix norms}\label{norms}
Recall the definitions \eqref{Alucia} and  \eqref{Adlucia} for $A(t,x)$ and $A_R(t,x)$ respectively. We work with the norm $\left\vert
y\right\vert _{A_R(t,x)}^{2}=\left\langle (A_RA^T_{R}(t,x))^{-1}y,y\right\rangle$,
$y\in \R^n$.

\label{page-norm1}

\begin{lemma}
\label{NORM1}
Let $x\in\R^n$, $t\geq 0$, $R>0$ and recall that $\lambda^*(A(t,x))$ and $\lambda_*(A(t,x))$ denote the largest and lowest singular value of $A(t,x)$.
\begin{itemize}
\item[$i)$] For every $y\in \R^n$ and $0<R\leq R^{\prime }\leq 1$
\begin{align}
\sqrt{\frac{R}{R^{\prime }}}\left\vert y\right\vert _{A_R(t,x)}&\geq \left\vert
y\right\vert _{A_{R^{\prime }}(t,x)}\geq \frac{R}{R^{\prime }}\left\vert
y\right\vert _{A_R(t,x)}\label{Norm2} \\
\frac{1}{\sqrt{R}
\lambda ^*(A(t,x))}\left\vert y\right\vert &\leq
\left\vert y\right\vert _{A_R(t,x)}\leq \frac{1}{R\lambda _*(A(t,x))}%
\left\vert y\right\vert .  \label{Norm3}
\end{align}

\item[$ii)$] For every $z\in \R^{m}$ and $R>0$
\begin{equation}
\left\vert A_R(t,x)z\right\vert _{A_R(t,x)}\leq \left\vert z\right\vert .
\label{Norm4}
\end{equation}

\item[$iii)$] For every $\varphi\in L^2([0,T];\R^m)$,
\begin{equation}
\Big|\int_0^r\varphi_s\,ds\Big|^2_{A_R(t,x)} \leq r\int_0^r|\varphi_s|^2_{A_R(t,x)}\,ds,\quad
r\in [0,T].  \label{Norm4bis}
\end{equation}
\end{itemize}
\end{lemma}

\bpr For fixed $x\in\R^n$ and $t\geq 0$, during the proof we omit in $A(t,x)$ and $A_R(t,x)$ the dependence on $(t,x)$, so we simply write $A$ and $A_R$

\smallskip

$i)$ For $0<R\leq R^{\prime }\leq 1$, it is easy to check that
\begin{equation*}
\frac{R^{\prime }}{R}A_RA_R^T\leq A_{R^{\prime }}A_{R^{\prime }}^T\leq \left( \frac{R^{\prime }}{R}\right) ^{2}A_RA_R^T
\end{equation*}%
which is equivalent to (\ref{Norm2}). This also implies (taking $
R^{\prime }=1$ so $A_{R^{\prime }}=A)$ that%
\begin{align*}
&\frac{1}{R}\lambda _*(A_R)^2\leq \lambda _*(A)^2\leq \frac{1}{R^{2}%
}\lambda _*(A_R)^2\\
&\frac{1}{R}\lambda ^*(A_R)^2 \leq \lambda ^*(A)^2\leq \frac{1}{R^{2}}\lambda ^*(A_R)^2
\end{align*}%
which immediately gives (\ref{Norm3}).

$ii)$ For $z\in \R^m$, we write $z=A_R^T y+w$ with $y\in \R^{n}$ and $%
w\in (\mathrm{Im}A_R^T)^{\bot }= \Ker A_R.$ Then $A_Rz=A_RA_R^T y$ so that%
\begin{eqnarray*}
\left\vert A_Rz\right\vert _{A_R}^{2} &=&\left\vert A_RA_R^T y\right\vert _{A_R}^{2}=\left\langle (A_RA_R^T)^{-1}A_RA_R^T y,A_RA_R^T y\right\rangle \\
&=&\left\langle y,A_RA_R^T y\right\rangle=\left\langle A_R^T y,A_R^T y\right\rangle =\left\vert A_R^T y\right\vert
^{2}\leq \left\vert z\right\vert ^{2}
\end{eqnarray*}
and (\ref{Norm4}) holds.

\smallskip

$iii)$ For $\varphi\in L^2([0,T];\R^m)$ and $r\in [0,T]$,
\begin{align*}
&\Big|\int_0^r\varphi_sds\Big|^2_{A_R} =\big\langle(A_RA_R^T)^{-1}\int_0^r\varphi_sds,%
\int_0^r\varphi_sds\big\rangle =\int_0^r\int_0^r\big\langle(A_RA_R^T)^{-1}\varphi_s,\varphi_u%
\big\rangle dsdu \\
&=\frac 12\int_0^r\int_0^r\big\langle(A_RA_R^T)^{-1}(\varphi_s-\varphi_u),\varphi_s-\varphi_u%
\big\rangle dsdu\\
&\quad-\int_0^r\int_0^r\Big(\big\langle(A_RA_R^T)^{-1}\varphi_s,\varphi_s\big\rangle -\big\langle(A_RA_R^T)%
^{-1}\varphi_u,\varphi_u\big\rangle \Big)dsdu \\
&=\frac 12\int_0^r\int_0^r\Big(|\varphi_s-\varphi_u|^2_{{A_R}}-2|\varphi_s|^2_{{A_R}}%
\Big)dsdu \\
&\leq  \int_0^r\int_0^r|\varphi_u|^2_{A_R}dsdu = t\int_0^t|\varphi_u|^2_{A_R}du.
\end{align*}
\epr

Next Lemma \ref{NORM3} is strictly connected to Remark \ref{1}, where we stressed that our result allows for a regime switch along the tube. In fact, here we fix $R>0$,  two points $(t,x)$ and $(s,y)$ and we get an equivalence between the norms $|\cdot|_{A_{R }(t,x)}$ and $|\cdot|_{A_{R }(s,y)}$ without supposing that in these two points the H\"ormander condition holds ``under the same regime''. To compensate this lack of uniformity, we suppose that the distance between $(t,x)$ and $(s,y)$ is bounded by $\sqrt{R}$, and we will need to take this fact into account. In the concatenation procedure of next Section \ref{sect-proof}, the size of the intervals, to which we apply our density estimates, will have to depend on the radius of the tube.

We set
$$
O=\{(t,x)\in [0,T]\times \R\,:\,\lambda(t,x)>0\}
$$
which is open, and under \eqref{Not2}, we define
\begin{equation}\label{DD}
\mathcal{D}=\Big\{C\,:\,O\to \R_+ \mbox{ such that } C(t,x)=K\Big(\frac{n(t,x)}{\lambda(t,x)}\Big)^q,\ K,q>0\Big\}.
\end{equation}
We also define
$$
1/\mathcal{D}=\Big\{c\,:\,O\to \R_+ \mbox{ such that } 1/c\in\mathcal{D}\}.
$$
\bl{NORM3}
Assume \eqref{Not2} and let $\mathcal{D}$ as in \eqref{DD}. There exists $C^*\in \mathcal{D}$ such that for every $(t,x),(s,y)\in O$ and $R\in(0,1]$ satisfying
\be{Norm5}
\left\vert x-y\right\vert +|t-s| \leq \sqrt{R}/C^*(t,x),
\ee
then for every $z\in \R^{n}$ one has
\be{Norm6}
\frac{1}{4}\left\vert z\right\vert^2 _{A_{R }(t,x)}\leq \left\vert
z\right\vert^2 _{A_{R }(s,y)}\leq 4\left\vert z\right\vert^2
_{A_{R }(t,x)}.
\ee
\el

\bpr \eqref{Norm6} is equivalent to
\[
4(A_{R }A_{R }^T)(t,x)\geq (A_{R }A_{R }^T)(s,y)\geq \frac{1}{4}(A_{R }A_{R }^T)(t,x),
\]
so we prove the above inequalities. Let $A_{R,k}$, $k=1,\ldots,m$, denote the columns of $A_R$. We use $(a+b)^{2}\geq \frac{1}{2}a^{2}-b^{2}$:
\begin{eqnarray*}
\langle A_{R }A_{R }^T(s,y)z,z\rangle
&=&\sum_{k=1}^{m}\left\langle A_{R ,k}(s,y),z\right\rangle ^{2} \\
&=&\sum_{k=1}^{m}(\left\langle A_{R ,k}(t,x),z\right\rangle
+\left\langle A_{R ,k}(s,y)-A_{R ,k}(t,x),z\right\rangle )^{2} \\
&\geq &\frac{1}{2}\sum_{k=1}^{m}\left\langle A_{R
,k}(t,x),z\right\rangle ^{2}-\sum_{k=1}^{m}\left\langle A_{R
,k}(s,y)-A_{R ,k}(t,x),z\right\rangle^{2}.
\end{eqnarray*}%
We use \eqref{Not2}: for every $(s,y)$ such that $|t-s|\leq 1$ and $|x-y|\leq 1$, we have
\begin{eqnarray*}
\langle A_{R }A_{R }^T(s,y)z,z\rangle
&\geq &\frac{1}{2}\sum_{k=1}^{m}\left\langle A_{R
,k}(t,x),z\right\rangle ^{2}-C_1 n(t,x)^\alpha R (\left\vert
x-y\right\vert ^{2}+\left\vert t-s\right\vert ^{2}) \left\vert
z\right\vert ^{2},
\end{eqnarray*}%
in which $C_1>0$ and $\alpha\geq 1$ denote universal constants. Notice that
\[
\sum_{k=1}^{m}\left\langle A_{R,k}(t,x),z\right\rangle ^{2}
=
\langle A_{R }A_{R }^T(t,x)z,z\rangle
\geq\l_*^2(A_{R}(t,x)) |z|^2
\geq R^2 \l_*^2(A(t,x))|z|^2.
\]
We choose the constants $(K,q)$ characterizing $C^*(t,x)$ such that $K\geq 2\sqrt{C_1}\vee 1$ and $q\geq \alpha$. So, under
\eqref{Norm5} we obtain
\[
C_1 n(t,x)^\alpha R (\left\vert
x-y\right\vert ^{2}+\left\vert t-s\right\vert ^{2}) \left\vert
z\right\vert ^{2} \leq \frac{1}{4}\sum_{k=1}^{m}\left\langle A_{R
,k}(t,x),z\right\rangle ^{2}
\]
and
\[
\left\langle (A_{R }A_{R }^T)(s,y)z,z\right\rangle \geq
\frac{1}{4}\sum_{k=1}^{m}\left\langle A_{R ,k}(t,x),z\right\rangle ^{2}=%
\frac{1}{4}\left\langle (A_{R }A_{R }^T)(t,x)z,z\right\rangle .
\]
The converse inequality follows from analogous computations and inequality $(a+b)^{2}\leq 2a^{2}+2b^{2}$.
\epr

We prove that moving along the skeleton associated to a control $\phi\in L^2([0,T],\R^d)$ for a small time $\delta$, the trajectory remains close to the initial point in the $A_\delta$-norm. To this purpose, we assume the conditions $(H_1)$ and $(H_2)$ in \eqref{Not6}. Notice that these give $(t,x_t(\phi))\in O$ for every $t$. Moreover, in such a case the set $\D$ can be replaced by the following class of functions:
\begin{equation}\label{AA}
\mathcal{A}=\Big\{C\,:\,[0,T]\to \R_+\,:\, C_t=K\Big(\frac{n_t}{\lambda_t}\Big)^q,\mbox{ for some $K,q>0$}\Big\},
\end{equation}
$n_t$ and $\lambda_t$ being defined in \eqref{Not6}. We also set
$$
1/\mathcal{A}=\big\{c\,:\,[0,T]\to (0,1]\,:\, 1/c_t\in \mathcal{A}\big\}.
$$

\bl{normcontrol}
Let $x(\phi)$ be the skeleton path \eqref{skeleton} associated to $\phi\in L^2([0,T],\R^d)$. Assume $(H_1)$ and $(H_2)$ in \eqref{Not6}.
Then there exists $\delta^*,\varepsilon^*\in 1/\mathcal{A}$ such that for every $t\in [0,T]$, $\delta_t\leq \delta^*_t$, $\varepsilon_t(\delta_t)\leq \varepsilon^*_t$, $ s\in[0, \delta_t]$ with $t+s\leq T$ and for every $z\in\R^n$ one has
\be{correction3}
\frac 14\left\vert z\right\vert^2 _{A_{\delta_t}(t,x_{t}(\phi))}\leq \left\vert
z\right\vert^2 _{A_{\delta_t}(t+s,x_{t+s}(\phi))}\leq 4\left\vert z\right\vert^2
_{A_{\delta_t}(t,x_{t}(\phi))}.
\ee

Moreover, there exists $\bar{C}\in \mathcal{A}$ such that
\be{normcontrolhat}
\sup_{0\leq s \leq \delta_t} |x_{t+s}(\phi)-(x_t(\phi)+b(t,x_t(\phi))s)|_{A_{\delta_t}(t,x_t(\phi))}\leq \bar{C}_t (\ve_t(\delta_t)\vee \sqrt{\delta_t})
\ee
where
\[
\ve_t(\delta)=\left(\int_t^{t+\delta} |\phi_s|^2ds\right)^{1/2}.
\]
\el

\bpr

Set $s_t=\inf\{s>0\,:\,|x_{t+s}(\phi)-x_t(\phi)|\geq 1\}$. From  \eqref{Not2} and $(H_1)$ in \eqref{Not6}, we have
\begin{align*}
1=|x_{t+s_t}(\phi)-x_t(\phi)|
&\leq n_t\big(s_t+\sqrt {s_t}\varepsilon_t(s_t)\big).
\end{align*}
We take $\underline{C}\in \mathcal{A}$ such that $n_t\big(\sqrt{s_t}+\varepsilon_t(s_t)\big)\leq \underline{C}_t^{1/2}$, so that $s_t\geq
1/\underline{C}_t$. Take now $\delta^*\in 1/\mathcal{A}$ such that $\delta^*\leq 1/\underline{C}$. Then if $s\leq \delta_t\leq \delta^*_t$, one has $s\leq s_t$ and
again from  \eqref{Not2} and $(H_1)$ in \eqref{Not6} we have
$$
|x_{t+s}(\phi)-x_t(\phi)|+|s|
\leq \sqrt{\delta_t}\big(n_t(\sqrt{\delta^*_t}\, +\varepsilon_t(\delta_t))+\sqrt{\delta^*_t}\big).
$$
By continuity, for every $\varepsilon^*\in1/\mathcal{A}$ and for every $t$ there exists $\hat\delta_t$ such that $\varepsilon_t(\hat \delta_t)\leq \varepsilon^*_t$. So, there actually exists  $\delta_t\leq \delta^*_t$ for which $\varepsilon_t(\delta_t)\leq \varepsilon^*_t$. And for such a $\delta_t$, we have
$$
|x_{t+s}(\phi)-x_t(\phi)|+|s|
\leq \sqrt{\delta_t}\big(n_t(\sqrt{\delta^*_t}\, +\varepsilon^*_t)+\sqrt{\delta^*_t}\big).
$$
We now choose $\delta^*,\varepsilon^*\in 1/\mathcal{A}$ in order that the last factor in the above right hand side is smaller than $1/C^*(t,x_t(\phi))$, where $C^*(t,x)$ is the function in $\mathcal{D}$ for which Lemma \ref{Norm3} holds.
Then \eqref{Norm5} is satisfied with $R=\delta_t$, $x=x_t(\phi)$, $y=x_{t+s}(\phi)$ and $s$ replaced by $t+s$. Hence \eqref{correction3} follows by applying \eqref{Norm6}.

\smallskip

We prove now \eqref{normcontrolhat}. For the sake of simplicity, we let $x_t$ denote the skeleton path $x_t(\phi)$. We write
\begin{align*}
J_{t,s}&:=x_{t+s}-x_{t}-b(t,x_{t})s=\int_{t}^{t+s}(\dot x_{u}-b(u,x_{u}))du +\int_{t}^{t+s}(b(u,x_{u})-b(t,x_{t}))du\\
&=\int_{t}^{t+s}\sigma(u,x_u)\phi_udu +\int_{t}^{t+s}(b(u,x_{u})-b(t,x_{t}))ds,
\end{align*}
so that
\begin{align*}
|J_{t,s}|^2_{A_{\delta_t}(t,x_t)}
&\leq 2s\int_t^{t+s}|\sigma(u,x_u)\phi_u|^2_{A_{\delta_t}(t,x_t)}dt
+2s\int_t^{t+s}|b(u,x_{u})-b(t,x_t)|^2_{A_{\delta_t}(t,x_t)}du.
\end{align*}
In the above right hand side, we apply \eqref{correction3} to the norm in the first term and we use \eqref{Norm3} in the second one. And we obtain:
\begin{align*}
|J_{t,s}|^2_{A_{\delta_t}(t,x_t)}
&\leq 2s\int_t^{t+s} 4|\sigma(u,x_u)\phi_u|^2_{A_{\delta_t}(u,x_{u})}du
+2s\int_t^{t+s}\frac{1}{\delta_t^2\lambda_t^2}|b(u,x_u)-b(t,x_t)|^2 du\\
&\leq 8s\int_t^{t+s} |\sigma(u,x_u)\phi_u|^2_{A_{\delta_t}(u,x_u)}du
+2\delta_t\int_t^{t+\delta_t}\frac{1}{\delta_t^2\lambda_t^2}\times n_t^2(|u-t|+ |x_u-x_t|)^2 du.
\end{align*}
We have already proved that, for $u\in[t,t+s]$, $|u-t|+ |x_u-x_t|\leq \sqrt{\delta_t}/C^*_t$, with $C^*\in \mathcal{A}$, so
\begin{align*}
|J_{t,s}|^2_{A_{\delta_t}(t,x_t)}
&\leq 8s\int_t^{t+s} |\sigma(u,x_u)\phi_u|^2_{A_\delta(u,x_u)}du
+\bar C_t\delta_t,
\end{align*}
with $\bar C\in \mathcal{A}$. It remains to study the first term in the above right hand side. For $i=1,\dots,m$, we set $\psi^{(j-1)d+j}=\frac{1}{\sqrt{\delta_t }}\phi^{j}$ for $j=1,\dots,d$, $\psi^i=0$ otherwise. Then, recalling \eqref{Adlucia}, we can write
$\s(u,x_{u})\phi_u=A_{\delta_t}(u,x_{u}(\phi)) \psi_u$,
so that, by \eqref{Norm4},
\[
|\sigma(u,x_u)\phi_u|^2_{A_{\delta_t}(u,x_u)}
=|A_{\delta_t}(u,x_u)\psi_u|^2_{A_{\delta_t}(u,x_u)}
\leq |\psi_u|^2=\frac 1{\delta_t}\,|\phi_u|^2.
\]
Hence, for $s\leq \delta_t$, we finally have
$|J_{t,s}|^2_{A_{\delta_t}(t,x_t)}\leq 8\ve_t(\delta_t)^2+\bar C_t\delta_t$, and the statement follows.
\epr

\begin{remark}\label{for-dc}
Let us finally discuss an inequality which will be used in next Section \ref{sectioncontrol}. Fix $x\in\R^n$ and let $x(\phi)$ be the skeleton path \eqref{skeleton} associated to $\phi\in L^2([0,T],\R^d)$ with starting condition $x_0(\phi)=x$. Assume simply \eqref{Not2} and recall $\mathcal{D}$ defined in \eqref{DD}. Then looking at the proof of Lemma \eqref{normcontrol}, we have the following result:
if $(0,x)\in O$, there exists $\overline{\delta}, \overline{\varepsilon}\in 1/\mathcal{D}$ and $\overline{C}\in\mathcal{D}$ such that if  $\delta\leq \overline{\delta}(0,x)$, $\varepsilon_0(\delta)\leq \overline{\varepsilon}(0,x)$ and $s\in[0, \delta]$ then
\be{normcontrolhat2}
\sup_{0\leq s \leq \delta} |x_{s}(\phi)-(x+b(0,x)s)|_{A_{\delta}(0,x)}\leq \overline{C}(0,x) (\ve_0(\delta)\vee \sqrt{\delta}).
\ee
\end{remark}

\subsection{Proof of Theorem \ref{mttubesS}}\label{sect-proof}

This section is organized as follows: the lower bound in Theorem \ref{mttubesS} is proved in next Theorem \ref{th-lower}, whereas the upper bound in Theorem \ref{mttubesS} is studied in next Theorem \ref{th-upper}.

As already mentioned, the proof we are going to develop relies on a two-sided bound for the density of equation \eqref{equation} in short time, proved in  \cite{BCP1}. The estimate is \emph{diagonal}, meaning that it is local around the drifted initial condition $x_0+b(0,x_0)\d$, $\d$ denoting the (small) time at which we are studying the density. But in order to be more precise and self-contained, we briefly recall the result from \cite{BCP1} we are going to strongly use.

We will suppose that
\begin{equation}\label{hpder}
\sum_{0\leq |\a| \leq 4} \Big[
\sum_{j=1}^d  |\partial_x^{\a} \s_j(t,x)|
+ |\partial_x^{\a} b(t,x)|
+ |\partial_x^{\a}\partial_t\s_j(t,x)| \Big]
\leq \kappa,
\ \forall t\in [0,T],\,\forall x \in \R^n.
\end{equation}
Of course \eqref{hpder} is much stronger than \eqref{Not2}  but we will see in the sequel that, by a suitable localization, one can reduce to the validity of \eqref{hpder} (see next Remark \ref{2}). We also assume that
\begin{equation}\label{ND}
\lambda (0,x_0)>0,
\end{equation}
$x_0$ denoting the starting point of the diffusion $X$ solving \eqref{equation}, and we consider the following set of constants:
\begin{equation}\label{CC}
\mathcal{D}_0=\Big\{C>0\,:\,C=K \Big(\frac{\kappa}{\lambda(0,x_0)}\Big)^{q},\ \exists\ K,q>0\Big\}.
\end{equation}
We use the notation $1/\mathcal{D}_0$ for constants $c$ such that $1/c\in\mathcal{D}_0$.

We set $p_\delta(x_0,\cdot)$ the density of $X_\delta$ when $X_0=x_0$. We use here the following version of Theorem \ref{m} in \cite{BCP1}:
\bt{m}
Suppose that \eqref{ND} and \eqref{hpder}  hold. Let $\D_0$ be defined in \eqref{CC}. Then there exist
$r^*,\d^*\in 1/\D_0$, $C\in \D_0$ such that for $\d\leq \d^*$ and for
$|y-x_0-b(0,x_0)\d|_{A_\d(0,x_0)} \leq r^*$  one has
\[
\frac{1}{ C \d^{n-\frac{\dim\langle \s(0,x_0) \rangle}{2} }}
\leq
p_\d(x_0,y)
\leq \frac{e^C}{  \d^{n-\frac{\dim\langle \s(0,x_0) \rangle}{2} }}
\]
where $\dim\langle \s(0,x_0) \rangle$ denotes the dimension of the vector space spanned by $\sigma_1(0,x_0),\ldots,$ $\sigma_d(0,x_0)$.
\et

Notice that \eqref{ND} and \eqref{hpder} are, respectively,   Assumption  2.2 and Assumption 2.3  in \cite{BCP1}. Therefore, Theorem \ref{m} is actually a re-writing of Theorem 3.7 in \cite{BCP1} (with the constant $C$ specified in Remark 3.8 therein) and Theorem 4.6 in \cite{BCP1}.

\begin{remark}\label{est-tx}
Of course, Theorem \ref{m} can be written for a general starting condition $(t,x)$ in place of $(0,x_0)$. In such a case, \eqref{ND} and \eqref{CC} have to be replaced by
$$
\lambda (t,x)>0\quad\mbox{and}\quad
\mathcal{D}_{t,x}=\Big\{C>0\,:\,C=K \Big(\frac{\kappa}{\lambda(t,x)}\Big)^{q},\ \exists\ c,q>0\Big\}
$$
respectively. But a closer look to the proofs of Theorem 3.7 and of Theorem 4.6 in \cite{BCP1} shows that the constants $K$ and $q$ in $\mathcal{D}_0$ are universal, that is, they can be taken independently of all the data (the starting point $(0,x_0)$, the diffusion coefficients, the quantities  $\lambda (0,x_0)$, $\kappa$ etc.). This means that Theorem \ref{m} can be formulated as follows. Assume that \eqref{hpder} holds and define the (open) set
$$
O=\big\{(t,x)\,:\,\lambda(t,x)>0\big\}.
$$
Set
$$
\mathcal{D}=\Big\{C:O\to\R_+\,:\,C(t,x)=K \Big(\frac{\kappa}{\lambda(t,x)}\Big)^{q},\ \exists\ c,q>0\Big\}.
$$
Then there exist $C\in \D$, $r^*,\d^*\in 1/\D$  such that for $(t,x)\in O$, $\d\leq \d^*(t,x)$ and for every $y$ such that
$|y-x-b(t,x)\delta|_{A_{\d}(t,x)} \leq r^*_t$  one has
$$
\frac{1}{ C(t,x) \d^{n-\frac{\dim\langle \s(t,x) \rangle}{2} }}
\leq
p(t,t+\delta,x,y)
\leq \frac{e^{C(t,x)}}{  \d^{n-\frac{\dim\langle \s(t,x) \rangle}{2} }},
$$
where $p(t,s,x,\cdot)$ denotes the density of the solution $X$ at time $s$ of the equation in \eqref{equation} but with the starting condition $X_t=x$. \end{remark}

\begin{remark}\label{C-B}
From \eqref{Adlucia} and the Cauchy-Binet formula we obtain (for details see (3.43) in  \cite{BCP1})
\be{CB}
\frac{1}{C(t,x)} \d^{n-\frac{\dim\langle \s(t,x) \rangle}{2}} \leq
\sqrt{\det A_\d A_\d^T (t,x)}\leq
C(t,x) \d^{n-\frac{\dim\langle \s(t,x) \rangle}{2}},
\ee
so the density bounds above  are  equivalent  to the following ones:
\be{denco}
\frac{1}{ C(t,x) \sqrt{\det A_{\d} A_{\d}^T (t,x)} }
\leq
p(t,t+\delta,x,y)
\leq \frac{e^{C(t,x)}}{  \sqrt{\det A_{\d} A_{\d}^T (t,x)}}
\ee
\end{remark}

\br{2}
The plan for the proof is the following. Consider first the lower bound (see Theorem \ref{th-lower}). For $\phi\in L^2[0,T]$, let $x(\phi)$ be the skeleton  associated to \eqref{equation} is given in \eqref{skeleton}. We set a discretization $0=t_0<t_1<\cdots <t_N=T$ of the time interval $[0,T]$. Then, as $k$ varies, we consider the events
\begin{equation}\label{D-Gamma}
D_k=\Big\{ \sup_{t_k\leq t \leq t_{k+1}} |X_t-x_t(\phi)|_{A_{R_t}(t,x_t(\phi))}\leq 1 \Big\} \mbox{ and }
\Gamma_k = \Big\{ y\,:\,|y-x_{t_k}(\phi)|_{A_{R_{t_k}}(t_k,x_{t_k}(\phi))} \leq r_k \Big\},
\end{equation}
where $r_k<1$ is a radius that will be suitably defined in the sequel. We denote $\PR_k$ the conditional probability
\[
\PR_k(\cdot)=\PR\left( \cdot|W_t,t\leq t_k; X_{t_k}\in \Gamma_k \right)
\]
We will lower bound
$\PR(\sup_{t\leq T} |X_t-x_t(\phi)|_{A_{R_t}(t,x_t(\phi))}\leq 1 )$ by
computing the product of the probabilities $\PR_k \left(D_k\cap \{X_{t_k+1}\in\Gamma_{k+1}\}\right) $, and this computation uses the lower estimate of the densities given in Theorem \ref{m}. Remark that Theorem \ref{m} uses \eqref{hpder}, a condition which asks for a global bound for the derivatives of the coefficients, whereas for the tube estimates we are assuming only $(H_1)$ in \eqref{Not6}, i.e. a bound for the coefficients which is not global but just in a neighborhood of the skeleton. Suppose that we have a process $X$ which, for some external reasons, verifies \eqref{equation} for $t_k\leq t \leq t_{k+1}$, and such that $\sup_{t_k\leq t \leq t_{k+1}} |X_t-x_t(\phi)|_{A_{R_t}(t,x_t(\phi))}\leq 1 $.
From $(H_1)$, $n_{t_k}$ bounds the derivatives of $\s(t,y)$ and $b(t,y)$ for all $(t_k-1)\vee 0\leq t \leq t_k+1$, and for all $|y-x_{t_k}(\phi)|\leq 1$.
Then, for example using the result in \cite{whitney1944}, we can define $\bar{\s},\bar{b}$ which coincide with $\s,b$ on $[(t_k-1)\vee 0, t_k+1]\times \{y\in \R^n: |y-x_{t_k}(\phi)|\leq 1\}$, are differentiable as many times as $\s,b$ but the bound in \eqref{hpder} holds on the whole $\R^+\times \R^n$. Let now $\bar{X}$ be the strong solution to
\[
\bar{X}_t=X_{t_k}+\int_{t_k}^t \bar{\s}(s,\bar{X}_s)\circ dW_s + \int_{t_k}^t \bar{b}(s,\bar{X}_s) ds, \quad t\in[t_k,t_{k+1}].
\]
Now, if we call $\bar D_k$ the sets in \eqref{D-Gamma} with $X$ replaced by $\bar X$, it is clear that
\[
\PR ( D_k  \cap \{X_{t_k+1}\in\Gamma_{k+1}\})=
\PR ( \bar D_k  \cap \{\bar{X}_{t_k+1}\in\Gamma_{k+1}\})
\]
and therefore we can equivalently prove our tube estimate supposing that the bound in $(H_1)$ holds globally, that is assuming \eqref{hpder}. This really allows us to apply Theorem \ref{m}. And a similar procedure can be developed for the upper bound (see Theorem \ref{th-upper}).
\er

We recall the set $\mathcal{A}$ defined in \eqref{AA}:
$$
\mathcal{A}=\Big\{C\,:\,[0,T]\to \R_+\,:\, C_t=K\Big(\frac{n_t}{\lambda_t}\Big)^q,\mbox{ for some $K, q>0$}\Big\}.
$$
We also recall $1/\mathcal{A}$ defined as usual. Notice that that, under \eqref{Not6}, $n(t,x_t(\phi))\leq n_t$ and $\lambda(t,x_t(\phi))\geq \lambda_t$. So, any $C(t,x)\in\mathcal{D}$ evaluated in $(t,x_t(\phi))$ is upper bounded by the function $C_t$ in $\mathcal{A}$ written with the same constants $K$ and $q$.

For $\mu\geq 1$, $h\in(0,1]$ and $K_*,q_*>0$, we denote
\be{Rmax}
R_t^*(\phi)=\exp\left(-K_* \left(\frac{\mu n_t}{\l_t}\right)^{q_*} \mu^{2q_*}\right)
\left( h\wedge \inf_{0\leq \d\leq h}
\left\{
\d \big/ \int_t^{t+\d} |\phi_s|^2 ds \right\}
\right)
\ee

\begin{theorem}\label{th-lower}
Let $\mu\geq 1$, $h\in (0,1]$, $n\,:\,[0,T]\to [1,+\infty)$, $\lambda\,:\,[0,T]\to (0,1]$, $\phi\in L^2([0,T],\R^n)$ and $R\,:\,[0,T]\to (0,1]$ be such that $(H_1)$--$(H_3)$ in \eqref{Not6} hold.
Then there exist $\bar K,\bar q>0$ such that
\begin{equation}\label{tube-lower}
\PR\left(\sup_{t\leq T} |X_t-x_t(\phi)|_{A_{R_t}(t,x_t(\phi))}\leq 1 \right)
\geq \exp\left(- \int_0^T \bar{K} \left(\frac{\mu n_t}{\l_t}\right)^{\bar{q}}\left( \frac 1h+\frac{1}{R_t}+|\phi_t|^2 dt\right) \right).
\end{equation}
Moreover, if $R_t \leq R_t ^*(\phi)$ for some $K_*,q_*>0$, $R^*(\phi)$ being given in \eqref{Rmax}, then
\begin{equation}\label{tube-lower-bis}
\PR\left(\sup_{t\leq T} |X_t-x_t(\phi)|_{A_{R_t}(t,x_t(\phi))}\leq 1 \right)
\geq \exp\left(- \int_0^T 2\bar{K} \left(\frac{\mu n_t}{\l_t}\right)^{\bar{q}}\left( \frac{1}{R_t}+|\phi_t|^2 dt\right) \right).
\end{equation}
\end{theorem}

\bpr
\noindent
\textbf{STEP 1.} We first set-up some quantities which will be used in the rest of the proof.

\smallskip

We recall $(H_{3})$: $R_.,\,|\phi _{.}|^{2},\,n_{.},\,\l_{.} \in L(\mu ,h)$, where
$f\in L(\mu ,h)$ if and only if $f(t)\leq \mu f(s)$ for $\left\vert t-s\right\vert \leq h$.
We set, for $q_1,K_1>1$ to be fixed in the sequel,
\[
f_R(t)= K_1 \left(\frac{\mu n_t}{\l_t}\right)^{q_1} \left( \frac{1}{h}+\frac{1}{R_t}+|\phi_t|^2\right).
\]
Then straightforward computations give that $f_R\in L(\mu^{2 q_1+1},h)$. We define
\be{tdelta}
\d(t)=\inf\left\{\d>0\,:\, \int_t^{t+\d} f_R(s) ds \geq \frac{1}{\mu^{2 q_1+1}} \right\}.
\ee
We have
\[
\frac{\d(t)}{h} = \int_t^{t+\d(t)} \frac{1}{h} ds \leq \int_t^{t+\d(t)} f_R(s)ds =\frac{1}{\mu^{2 q_1+1}},
\]
so $\d(t)\leq h$.
%, and we can use on the intervals $[t,t+\d(t)]$ the fact that the functions $n_\cdot,\lambda_\cdot, R_\cdot, |\phi_\cdot|^2$ are in $L(\mu,h)$.
We now prove that $\delta(\cdot)\in L(\mu^{4q_1+1}, h)$. In fact, if $0<t-t'\leq h$,
\[
\mu^{2q_1+1} f_R(t) \d(t) \geq \int_{t}^{t+\delta(t)} f_R (s) ds=\frac 1{\mu^{2q_1+1}}
=\int_{t'}^{t'+\d(t')}f_{R}(s)ds \geq \mu^{-(2q_1+1)} f_{R}(t)\d(t'),
\]
so $\d(t')\leq \mu^{4q_1+2}\d(t)$. Since the converse holds as well, we get
$\d(\cdot) \in L(\mu^{4q_1+2},h)$. We now prove a further property for $\delta(\cdot)$:
we have
\[
\frac{1}{\mu^{2 q_1+1}}= \int_t^{t+\d(t)} f_R(s) ds
\geq \int_t^{t+\d(t)} \frac{f_R(t)}{\mu^{2 q_1+1}} ds
\geq \d(t) \frac{f_R(t)}{\mu^{2 q_1+1}},
\]
so
\be{Rdelta}
\d(t) \leq \frac{1}{f_R(t)} \leq \frac{R_t}{K_1} \left(\frac{\l_t}{\mu n_t}\right)^{q_1}
\leq \frac{1}{K_1} \left(\frac{\l_t}{\mu n_t}\right)^{q_1} \in1/\mathcal{A}
\ee
(recall that $R_t,\lambda_t\leq 1$ and $n_t\geq 1$ for every $t$). We also set the energy over the time interval $[t,t+\delta(t)]$:
\[
\ve_t(\delta(t))=\left( \int_t^{t+\delta(t)} |\phi_s|^2 ds\right)^{1/2}.
\]
Since $n,\lambda\in L(\mu, h)$ and $\delta(t)\leq h$, for $s\in (t,t+\delta(t))$ we have
$$
f_R(s)\geq K_1\Big(\frac{\mu n_s}{\lambda_s}\Big)^{q_1}|\phi_s|^2
\geq \frac{K_1}{\mu^{2q_1}}\Big(\frac{\mu n_t}{\lambda_t}\Big)^{q_1}|\phi_s|^2.
$$
Hence
$$
\frac 1{\mu^{2q_1+1}}=\int_t^{t+\delta(t)}f_R(s)ds\geq
\frac{K_1}{\mu^{2q_1}}\Big(\frac{\mu n_t}{\lambda_t}\Big)^{q_1}\int_t^{t+\delta(t)}|\phi_s|^2ds,
$$
which gives that
\be{energy}
\ve_t(\delta(t))^2 \leq \frac{1}{K_1} \left(\frac{\l_t}{\mu n_t}\right)^{q_1}\in1/\mathcal{A}.
\ee

\textbf{STEP 2.}
We set now some notation and properties that will be used in the ``concatenation'', which is developed in the following steps.

\smallskip

We define the time grid as
\[
t_0=0,\quad t_k=t_{k-1}+\d(t_{k-1}),
\]
and introduce the following notation on the grid:
\[
\d_k=\d(t_k),\ \  \ve_k=\ve_{t_k}(\delta_k),\ \  n_k=n_{t_k},\ \  \lambda_k= \lambda_{t_k},\ \  X_k=X_{t_k},\ \  x_k=x_{t_k}(\phi),\ \  R_k=R_{t_k}.
\]
Recall that $\delta(t)<h$ for every $t$, so we have
\[
R_k/\mu\leq R_t\leq \mu R_k,\quad \mbox{ for } t_k\leq t \leq t_{k+1}.
\]
We also define
\[
\hat{X}_k= X_k + b(t_k,X_k)\delta_k,
\quad
\hat{x}_k=x_k+b(t_k, x_k)\delta_k,
\]
and for $t_k\leq t \leq t_{k+1}$,
\[
\hat{X}_k(t)=X_k+b(t_k,X_k)(t-t_k),
\quad
\hat{x}_k(t)=x_k+b(t_k, x_k)(t-t_k).
\]
Let $r^*\in1/\mathcal{A}$ be the radius-function of Theorem \ref{m}, in the version of Remark \ref{est-tx}, associated to the points $(t,x_t(\phi))$ as $t\in[0,T]$. We set $r^*_k=r^*_{t_k}$.
%
%Moreover we denote $r^*_k$
%\[
%\cal{C}_k=\cal{C}_{t_k},
%\]
%and $r_k^*\in \cal{C}_k$ the radius $r^*$ of Theorem \ref{m} associated to the diffusion starting at time $t_k$ (instead of at time 0).
%

Let us see some properties.

For all $t_k\leq t\leq t_{k+1}$, we have $R_t\geq R_k/\mu\geq \delta_k/\mu$ and, by using \eqref{Norm2}, we obtain
$$
|\xi|_{A_{R_t}(t,x_t)} \leq \sqrt{\frac{\delta_k}{R_k}}\,
|\xi|_{A_{\delta_k/\mu}(t,x_t)}
\leq |\xi|_{A_{\delta_k/\mu}(t,x_t)},
$$
last inequality holding because $\delta_k\leq R_k$. Since $\delta_k/\mu\leq \delta_k$, we apply again \eqref{Norm2} to the norm in the right hand side above and we get
\begin{equation}\label{magia0}
|\xi|_{A_{R_t}(t,x_t)} \leq \mu\,
|\xi|_{A_{\delta_k}(t,x_t)}.
\end{equation}
Taking $\xi=x_t-\hat{x}_k(t)$, we have
$$
|x_t-\hat{x}_k(t)|_{A_{R_t}(t,x_t)} \leq \mu\,
|x_t-\hat{x}_k(t)|_{A_{\delta_k}(t,x_t)}.
$$
By  \eqref{Rdelta} and \eqref{energy}, we can choose $q_1, K_1$ large enough
such that $\delta(t)\leq \delta^*(t)$, $\ve_t(\delta(t))\leq \ve^*(t)$ where $\delta^*\in 1/\mathcal{A}$ and $\ve^*\in 1/\mathcal{A}$ are the functions in Lemma \ref{normcontrol}. So, we apply \eqref{correction3} to the norm in the above right hand side and we obtain
$$
|x_t-\hat{x}_k(t)|_{A_{R_t}(t,x_t)} \leq \mu
\times 4|x_t-\hat{x}_k(t)|_{A_{\delta_k}(t_k,x_k)}.
$$
We use now \eqref{normcontrolhat}: for some $\bar C\in \mathcal{A}$, we get
$$
|x_t-\hat{x}_k(t)|_{A_{\delta_k}(t_k,x_k)}
\leq \bar C_k(\ve_k\vee\sqrt{\d_k})
$$
where $\bar C_k=\bar C_{t_k}$, and, as a consequence of the estimate above, we have also
$$
|x_t-\hat{x}_k(t)|_{A_{R_t}(t,x_t)} \leq 4\mu \,\bar C_k(\ve_k\vee\sqrt{\d_k}),
$$
for all $t\in[t_k,t_{k+1}]$ and for all $k$. By recalling that
$x_{t_{k+1}}-\hat{x}_k(t_{k+1})=x_{k+1}-\hat{x}_k$, and possibly choosing $K_1$ larger, we can resume by asserting that $\d_k\leq \d^*_{t_k}$ in Theorem \ref{m} with initial condition $(t_k,x_k)$ (see its version in Remark \ref{est-tx})  and
\begin{align}
|x_{k+1}-\hat{x}_k|_{A_{\delta_k}(t_k,x_k)} &\leq r_k^*/4 \mbox{ for all } k,\label{magia1}\\
|\hat{x}_k(t)-x_t|_{A_{R_t}(t,x_t)} &\leq \frac{1}{4} \mbox{ for all } t\in[t_k, t_{k+1}] \mbox{ and for all } k.\label{magia2}
\end{align}
We have already noticed that, under our settings, \eqref{correction3} holds, so that
\[
\frac{1}{2} |\xi|_{A_{\delta_k}(t_k,x_k)}\leq |
\xi|_{A_{\delta_k}(t_{k+1},x_{k+1})}\leq 2 |\xi|_{A_{\delta_k}(t_k,x_k)}.
\]
Since $\d(\cdot) \in L(\mu^{4q_1+2},h)$, one has $\d_k/\d_{k+1}\leq
\mu^{4q_1+2}$ and $\d_{k+1}/\d_k\leq \mu^{4q_1+2}$. So, using \eqref{Norm2} to the right hand side of the above inequality we easily get
\be{fe}
\frac{1}{2\mu^{2q_1+1}}|\xi|_{A_{\delta_k}(t_k,x_k)} \leq |\xi|_{A_{\delta_{k+1}}(t_{k+1},x_{k+1})}\leq 2\mu^{2q_1+1} |\xi|_{A_{\delta_k}(t_k,x_k)}
\mbox{ for all }k.
\ee

\textbf{STEP 3.} We are ready to set-up the concatenation for the lower bound.

\smallskip

We set, for $K_2$ and $q_2$ to be fixed in the sequel,
\be{ray}
r_k= \frac{1}{K_2 \mu^{2q_1+2q_2+1} }\left(\frac{\l_k}{n_k}\right)^{q_2}.
\ee
Moreover, since $\lambda,n\in L(\mu,h)$ and $\delta_k\leq h$, one easily gets $r_{k+1}/r_k\leq \mu^{2 q_2}$ for every $k$.

%We define
%\[
%\Gamma_k=\{ |X_k-x_k|_{A_{\delta_k}(t_k,x_k)} \leq r_k \},\quad D_k=\{ \sup_{t_k\leq t \leq t_{k+1}} |X_t-x_t|_{A_{R_t}(t,x_t)}\leq 1 \},
%\]
%and $\PR_k$ as the conditional probability
%\[
%\PR_k(\cdot)=\PR\left( \cdot|W_t,t\leq t_k; X_k\in \Gamma_k \right).
%\]

We define
\[
\Gamma_k=\{y\,:\, |y-x_k|_{A_{\delta_k}(t_k,x_k)} \leq r_k \}\mbox{ and }
\PR_k(\cdot)=\PR\left( \cdot|W_t,t\leq t_k;  X_k\in\Gamma_k \right),
\]
that is, $\PR_k$ is the conditional probability with respect to the knowledge of the Brownian motion up to time $t_k$ and the fact that $X_k\in\Gamma_k$.
The aim of this step is to prove that
\begin{equation}\label{to-prove1}
\PR_k(X_{k+1}\in\Gamma_{k+1})\geq 2 \mu^{-4n q_1} \exp(-K_3 (\log \mu + \log n_k -\log \l_k ))\mbox{ for all }k.
\end{equation}
for some constant $K_3$ depending on $K_1$, $K_2$, $q_1$ and $q_2$.

We denote $\rho_k(X_k,y)$ the density of $X_{k+1}$ with respect to this probability.
We prove that
\begin{equation}\label{to-prove2}
\Gamma_{k+1}
%\{y: \, |y -x_{k+1}|_{{A_{\delta_{k+1}}(t_{k+1},x_{k+1})}}\leq r_{k+1}\}
\subset\{y:\,|y-\hat X_k|_{A_{\delta_k}(t_k,X_k)}\leq r^*_k\}.
\end{equation}
If \eqref{to-prove2} holds, as we will see, then we can apply the lower bound in Remark \ref{est-tx} to $\rho_k(X_k,y)$. More precisely, we use here the version of the estimate given in \eqref{denco}: there exists $\underline{C}\in\mathcal{A}$ such that
\be{eqsttt}
\rho_k(X_k,y)\geq \frac{1}{\underline{C}_k \sqrt{\det A_{\d_k} A_{\d_k}^T(t_k,X_k)}}\mbox{ for all }y\in \Gamma_{k+1},
\ee
where $\underline{C}_k=\underline{C}_{t_k}$. Let us show that \eqref{to-prove2} holds. We estimate
$$
|y-\hat{X}_{k}|_{A_{\delta_{k}}(t_{k},x_{k})} \leq
|y-x_{k+1}|_{A_{\delta_{k}}(t_{k},x_{k})}+|x_{k+1}-\hat{x}_k|_{A_{\delta_{k}}(t_{k},x_{k})}
+|\hat{x}_k-\hat{X}_k|_{A_{\delta_{k}}(t_{k},x_{k})}
$$
and by using \eqref{magia1} we obtain
\be{deco}
|y-\hat{X}_{k}|_{A_{\delta_{k}}(t_{k},x_{k})} \leq
|y-x_{k+1}|_{A_{\delta_{k}}(t_{k},x_{k})}+\frac{r^*_k}4
+|\hat{x}_k-\hat{X}_k|_{A_{\delta_{k}}(t_{k},x_{k})}.
\ee
Using \eqref{fe}, the fact that $r_{k+1}/r_k\leq \mu^{2 q_2}$ and recalling that $|y -x_{k+1}|_{A_{\delta_{k+1}}(t_{k+1},x_{k+1})}\leq r_{k+1}$, we obtain
\begin{align*}
|y-x_{k+1}|_{A_{\delta_{k}}(t_{k},x_{k})}
&\leq 2 \mu^{2 q_1+1} |y-x_{k+1}|_{A_{\delta_{k+1}}(t_{k+1},x_{k+1})}
\leq 2 \mu^{2 q_1+1} r_{k+1}\\
&\leq 2 \mu^{2 q_1+2 q_2 + 1} r_{k}
\leq  \frac{2}{K_2}\left(\frac{\l_k}{n_k}\right)^{q_2}.
\end{align*}
\eqref{Not2}  also gives $|\hat{x}_k-\hat{X}_k|_{A_{\delta_{k}}(t_{k},x_{k})} \leq C_k |{x}_k-{X}_k|_{A_{\delta_{k}}(t_{k},x_{k})}$, where $C_k=C_{t_k}$ and $C$ is a suitable function in $\mathcal{A}$, and the conditioning with respect to $\Gamma_k$ gives $|\hat{x}_k-\hat{X}_k|_{A_{\delta_{k}}(t_{k},x_{k})} \leq C_k r_k$. Similarly,
$|\hat{x}_k(t)-\hat{X}_k(t)|_{A_{R_t}(t,x_t)} \leq C_k|x_k-X_k|_{A_{R_t}(t,x_t)}$ and by using firstly \eqref{magia0} and secondly \eqref{correction3}, we get
$$
|\hat{x}_k(t)-\hat{X}_k(t)|_{A_{R_t}(t,x_t)} \leq C_k\times \mu
|x_k-X_k|_{A_{\delta_k}(t,x_t)}
\leq C_k\mu\times 2
|x_k-X_k|_{A_{\delta_k}(t_k,x_k)}
\leq 2\mu C_k r_k,
$$
for every $t\in[t_k,t_{k+1}]$. Recalling \eqref{ray}, $K_2$ and $q_2$ (possibly large) such that
$|y-x_{k+1}|_{A_{\delta_{k}}(t_{k},x_{k})}\leq r^*_k/8$, $|\hat{x}_k-\hat{X}_k|_{A_{\delta_{k}}(t_{k},x_{k})} \leq r^*_k/8$, and
\be{magia3}
|\hat{X}_k(t)-\hat{x}_k(t)|_{A_{R_t}(t,x_t)} \leq 1/4,\quad \mbox{ for all } t\in[t_k, t_{k+1}]\mbox{ and for all }k.
\ee
From \eqref{deco}, this implies $|y-\hat{X}_{k}|_{A_{\delta_{k}}(t_{k},x_{k})} \leq r^*_k/2$.
On the event $\Gamma_k$, we also have, from \eqref{Norm3}, $|x_k-X_k|\leq |x_k-X_k|_{A_{\delta_{k}}(t_{k},x_{k})} \lambda^*(A(t_k,x_k))\sqrt{\d_k}\leq n_{t_k}^\alpha \sqrt{\d_k}\, r_k$, for some universal constant $\alpha>0$. So, we can fix $K_2$ and $q_2$ in order that Lemma \ref{NORM3} holds with $R=\delta_k$, $x=x_k$, $y=X_k$, $t=t_k$ and $s=0$. Then, we get
\[
\frac{1}{2}|\xi|_{A_{\delta_{k}}(t_{k},x_{k})}\leq|\xi|_{A_{\d_k}(t_k,X_k)}\leq 2|\xi|_{A_{\delta_{k}}(t_{k},x_{k})}.
\]
These inequalities give two consequences. First, we have
$$
|y-\hat{X}_{k}|_{A_{\d_k}(t_k,X_k)} \leq 2
|y-\hat{X}_{k}|_{A_{\d_k}(t_k,x_k)}\leq  r^*_k,
$$
so that \eqref{to-prove2} actually holds and then \eqref{eqsttt} holds as well. As a second consequence, we have that
\begin{align*}
&\Big\{y:|y -x_{k+1}|_{A_{\d_k}(t_k,X_k)}\leq \frac{r_{k+1}}{4\mu^{2 q_1+1}} \Big\}
\subset \Big\{y:|y -x_{k+1}|_{A_{\delta_{k}}(t_{k},x_{k})}\leq \frac{r_{k+1}}{2\mu^{2 q_1+1} }\Big\}\\
&\quad\subset \{y:|y -x_{k+1}|_{A_{\delta_{k+1}}(t_{k+1},x_{k+1})}\leq r_{k+1}\} =\Gamma_{k+1}  ,
\end{align*}
in which we have used \eqref{fe}.
Since $ r_{k+1}/(4\mu^{2 q_1+1} )\geq r_k /(4\mu^{2 q_1+ 2q_2+1} )$, we obtain
\begin{align*}
\Gamma_{k+1}
\supset
\Big\{y:|y -x_{k+1}|_{A_{\d_k}(t_k,X_k)}\leq \frac{r_k}{ 4\mu^{2 q_1+ 2q_2+1}}
\Big\}.
\end{align*}
By recalling that $r_k /(4\mu^{2 q_1+ 2q_2+1} )= \frac{1}{4 K_2 \mu^{4 q_1+ 4 q_2+2} }\left(\frac{\l_k}{n_k}\right)^{q_2}$, we can write, with $\Leb_n$ denoting the Lebesgue measure in $\R^n$,
\[
\Leb_n(\Gamma_{k+1}) \geq
\sqrt{\det(A_{\d_k} A_{\d_k}^T(t_k,X_k))}
\left(\frac{1}{4 K_2 \mu^{4 q_1+ 4 q_2+2} }\left(\frac{\l_k}{n_k}\right)^{q_2}\right)^n.
\]
So, from \eqref{eqsttt},
\[
\PR_k(X_{k+1}\in \Gamma_{k+1})\geq \frac{1}{\underline{C}_k} \left(\frac{1}{4 K_2 \mu^{4 q_1+ 4 q_2+2} }\left(\frac{\l_k}{n_k}\right)^{q_2}\right)^n
\]
where $\underline{C}_k$ is the constant in \eqref{eqsttt}. This implies \eqref{to-prove1}, for some constant $K_3$ depending on $K_2$ and $q_2$.

\smallskip

\textbf{STEP 4.} We give here the proof of the lower bounds \eqref{tube-lower} and
\eqref{tube-lower-bis}.

We set
\[
D_k=\Big\{ \sup_{t_k\leq t \leq t_{k+1}} |X_t-x_t|_{A_{R_t}(t,x_t)}\leq 1 \Big\}
\mbox{ and }
E_k=\Big\{\sup_{t_k\leq t\leq t_{k+1}} |X_t-\hat{X}_k(t)|_{A_{R_t}(t,x_t)}\leq \frac{1}{2}\Big\}.
\]
For $t\in[ t_k, t_{k+1}]$, by using \eqref{magia2} and \eqref{magia3} we have
\begin{align*}
|X_t-x_t|_{A_{R_t}(t,x_t)}
&\leq|X_t-\hat{X}_k(t)|_{A_{R_t}(t,x_t)}+|\hat{X}_k(t)-\hat{x}_k(t)|_{A_{R_t}(t,x_t)}
+|\hat{x}_k(t)-x_t|_{A_{R_t}(t,x_t)}\\
&\leq |X_t-\hat{X}_k(t)|_{A_{R_t}(t,x_t)}+\frac 12,
\end{align*}
so that
$E_k\subset D_k$. Moreover, by passing from Stratonovich to It\^o integrals and by using \eqref{Norm3}, we have
\[
\begin{split}
&|X_t-\hat{X}_k(t)|_{A_{R_t}(t,x_t)}
\leq |\s(t_k,X_{t_k})(W_t-W_{t_k})|_{A_{R_t}(t,x_t)}\\
&\quad + \Big|\int_{t_k}^{t}\big(\s(s,X_s)-\s(t_k,X_k)\big) dW_s\Big|_{A_{R_t}(t,x_t)}+ \Big| \int_{t_k}^{t} \big(b(s,X_s)-b(t_k,X_k)\big) ds\Big|_{A_{R_t}(t,x_t)}\\
&\quad+ \sum_{l=1}^d \Big|\int_{t_k}^{t} \nabla \s_l(s,X_s)(\s_l(s,X_s)-\s_l(t_k,X_k)) ds\Big|_{A_{R_t}(t,x_t)}\\
&\quad\leq  \Big| \frac{\sqrt{\mu}}{\sqrt{R_k}}\s(t_k,X_{t_k})(W_t-W_{t_k})\Big|_{A(t,x_t)}+ \Big|\frac{\mu}{R_k}\int_{t_k}^{t}\big(\s(s,X_s)-\s(t_k,X_k)\big) dW_s\Big|\\
&\quad+  \Big|\frac{\mu}{R_k} \int_{t_k}^{t} \big(b(s,X_s)-b(t_k,X_k)\big) ds\Big|+ \sum_{l=1}^d \Big|\frac{\mu}{R_k}\int_{t_k}^{t} \frac{\nabla \s_l(s,X_s)}{2}(\s_l(s,X_s)-\s_l(t_k,X_k)) ds\Big|.
\end{split}
\]
We use now the exponential martingale inequality (see also Remark \ref{2}) and we find that
\[
\PR_k(E_k^c) \leq  \exp\left(- \frac{1}{K_4} \left(\frac{\l_k}{\mu n_k}\right)^{q_4} \frac{R_k}{\d_k} \right)
\]
for some constants $K_4,q_4$. From \eqref{Rdelta}, $R_k/\d_k\geq K_1(\mu n_k/\l_k)^{q_1}$, so by choosing $K_1$ and $q_1$ possibly larger  and by recalling \eqref{to-prove1}, we can conclude that
\[
\PR_k(E^c_k) \leq \mu^{-4 nq_1} \exp(-K_3 (\log \mu + \log n_k -\log \l_k )) \leq \frac{1}{2}\PR_k(X_{k+1}\in\Gamma_{k+1}).
\]
Hence,
\be{shorttimediff}
\begin{split}
&\PR_k(\{X_{k+1}\in\Gamma_{k+1}\}\cap D_k) \geq \PR_k(\{X_{k+1}\in\Gamma_{k+1}\}\cap E_k)\geq \PR_k(X_{k+1}\in\Gamma_{k+1})-\PR_k(E^c_k)\\
&\quad \geq \frac{1}{2} \PR_k(X_{k+1}\in\Gamma_{k+1})\geq \exp\left( -K_5 (\log \mu +\log n_k -\log \l_k) \right),
\end{split}
\ee
for some constant $K_5$.
Let now $N(T)=\max\{k: t_k\leq T\}$. From definition \eqref{tdelta},
\[
\int_0^T f_R(t)dt \geq \sum_{k=1}^{N(T)} \int_{t_{k-1}}^{t_k} f_R(t)dt =
\frac{N(T)}{\mu^{2q_1+1}}.
\]
From \eqref{shorttimediff},
\begin{align*}
&\PR\Big(\sup_{t\leq T} |X_t-x_t|_{A_{R_t}(t,x_t)}\leq 1 \Big)
 \geq \PR\Big(\bigcap_{k=1}^{N(T)} \{X_{k+1}\in\Gamma_{k+1}\}\cap D_k\Big) \\
&\quad  \geq \prod_{k=1}^{N(T)} \exp(-K_5(\log \mu +\log n_k -\log \l_k))
= \exp\Big(-K_5\sum_{k=1}^{N(T)}\big(\log \mu +\log n_k -\log \l_k\big)\Big).
\end{align*}
Since
\[
\begin{split}
\sum_{k=1}^{N(T)}(\log \mu +\log n_k -\log \l_k)
&=\mu^{2q_1+1}\sum_{k=1}^{N(T)}\int_{t_k}^{t_k+1} f_R(t)(\log \mu +\log n_k -\log \l_k)dt \\
&\leq \mu^{2q_1+1} \int_0^T f_R(t) \log\Big(\frac{\mu^{3} n_t}{\l_t}\Big)dt,
\end{split}
\]
the lower bound \eqref{tube-lower} follows. Concerning \eqref{tube-lower-bis}, it is immediate from \eqref{tube-lower} and the fact that $R_t \leq R_t ^*(\phi)\leq h\, \exp\left(-K_* \left(\frac{\mu n_t}{\l_t}\right)^{q_*}\right)$.
\epr

We can now address the problem of the upper bound.

\begin{theorem}\label{th-upper}
Let $\mu\geq 1$, $h\in (0,1]$, $n\,:\,[0,T]\to [1,+\infty)$, $\lambda\,:\,[0,T]\to (0,1]$, $\phi\in L^2([0,T],\R^n)$ and $R\,:\,[0,T]\to (0,1]$ be such that $(H_1)$--$(H_3)$ in \eqref{Not6} hold. Suppose that, for some $K_*,q_*>0$ and for $R^*(\phi)$ as in \eqref{Rmax}, one has $R_t \leq R_t ^*(\phi)$. Then there exist $\bar K,\bar q>0$ such that \begin{equation}\label{tube-upper}
\begin{split}
&\PR\Big(\sup_{t\leq T} |X_t-x_t(\phi)|_{A_{R_t}(t,x_t(\phi))}\leq 1 \Big)\\
&\qquad\leq \exp\Big(- \int_0^T \bar{K} \Big(\frac{\mu n_t}{\l_t}\Big)^{\bar{q}}
\Big[\frac{\exp\Big(-K_* \Big(\frac{\mu n_t}{\l_t}\Big)^{q_*}\Big)}{R_t}+|\phi_t|^2\Big]dt \Big).
\end{split}
\end{equation}
\end{theorem}
\bpr
We refer here to notation and arguments already introduced and developed in the proof of Theorem \ref{th-lower}. So, when we recall here STEP 1, 2 and 3, we intend to refer to the same steps developed in the proof of Theorem \ref{th-lower}.

We define, with the same $K_1,q_1$ as in STEP 1,
\[
g_R(t)=
K_1 \left(\frac{\mu n_t}{\l_t}\right)^{q_1}
\left( \frac{\exp \left( -K_*\left(\frac{\mu n_t}{\l_t}\right)^{q_*} \mu^{2q_*} \right)}{R_t}
+|\phi_t|^2 \right)
\]
Because of \eqref{Rmax}, for all $t\in[0,T$],
\be{rh}
\frac{\exp \left( -K_*\left(\frac{\mu n_t}{\l_t}\right)^{q_*}\mu^{2q_*} \right)}{R_t} \geq \frac{1}{h}
\ee
We work here with  $\d(t)$ as in the proof od Theorem \ref{th-lower} but defined from $g_R$: \[
\d(t)= \inf\Big\{\d>0\,:\, \int_t^{t+\d} g_R(s) ds \geq \frac{1}{\mu^{2q_1+1}} \Big\}.
\]
We set, as before,
\[
\ve_t(\delta(t))=\Big( \int_t^{t+\delta(t)} |\phi_s|^2 ds\Big)^{1/2}.
\]
As in STEP 1, using also \eqref{rh}, we can check estimates similar to \eqref{Rdelta} and \eqref{energy}: we have indeed,
\[
\d(t) \leq \frac{h}{K_1} \Big(\frac{\l_t}{\mu n_t}\Big)^{q_1} \leq \frac{1}{K_1} \Big(\frac{\l_t}{\mu n_t}\Big)^{q_1}\quad\mbox{ and }\quad
\ve_t(\delta(t))^2 \leq \frac{1}{K_1} \Big(\frac{\l_t}{\mu n_t}\Big)^{q_1}.
\]
In particular, $\d(t)\leq h$. With these definitions we set a time grid $\{t_k: k=0,\dots, N(T)\}$ and all the associated quantities as in STEP 2.  As we did for the lower bound, since we estimate the probability of remaining in the tube for any $t\in [t_k,t_{k+1}]$, we can suppose that the bound in \eqref{hpder} holds on $ \R^+\times \R^n$ (recall Remark \ref{2}).
The short time density estimate \eqref{denco} holds again.
Recall now that $R_{.}\in L(\mu,h)$, and this gives the analogous to \eqref{fe}:
\be{fewR}
\frac{1}{2 \sqrt{\mu}}|\xi|_{A_{R_k}(t_k,x_k)}\leq
|\xi|_{A_{R_{k+1}}(t_{k+1},x_{k+1})} \leq
2 \sqrt{\mu} |\xi|_{A_{R_k}(t_k,x_k)}
\ee
We define
\[
\Delta_k=\{y\,:\, |y-x_k|_{A_{R_k}(t_k,x_k)} \leq 1 \}
\mbox{ and } \tilde{\PR}_k(\cdot)=\PR\left( \cdot|W_t,t\leq t_k; X_k\in \Delta_k \right),
\]
so $\tilde{\PR}_k$ is the conditional probability given the Brownian path up to time $t_k$ and the fact that $X_k\in \Delta_k$.

Now, since $\d(t)\leq h$ and $R,\l,n\in L(\mu,h)$, we have
$$
\int_t^{t+\d(t)}  K_1 \left(\frac{\mu n_s}{\l_s}\right)^{q_1} |\phi|_s^2 ds \leq
\mu^{2 q_1} K_1 \left(\frac{\mu n_t}{\l_t}\right)^{q_1} \int_t^{t+\d(t)}   |\phi|_s^2 ds
$$
and
\[
\begin{split}
&\int_t^{t+\d(t)}  K_1 \left(\frac{\mu n_s}{\l_s}\right)^{q_1}
\frac{\exp \left( -K_*\left(\frac{\mu n_s}{\l_s}\right)^{q_*} \mu^{2q_*}\right) }{R_s} ds\\
&\qquad\qquad\qquad\leq
\mu^{2q_1+1} K_1 \left(\frac{\mu n_t}{\l_t}\right)^{q_1}\exp \left( -K_*\left(\frac{\mu n_t}{\l_t}\right)^{q_*} \right) \frac{\d(t)}{R_t}.
\end{split}
\]
Since
\[
R_t\leq R_t^*(\phi)\leq\exp\left(-K_* \left(\frac{\mu n_t}{\l_t}\right)^{q_*} \mu^{2q_*}\right)
\left( \inf_{0\leq \d\leq h}
\left\{
\d \big/ \int_t^{t+\d} |\phi_s|^2 ds \right\}
\right),
\]
we have
\[
\int_t^{t+\d(t)} |\phi_s|^2 ds\leq
\exp\left(-K_* \left(\frac{\mu n_t}{\l_t}\right)^{q_*} \right)\frac{\d(t)}{R_t}
\]
We obtain
\[
1=
\mu^{2q_1+1}\int_t^{t+\d(t)}  g_R(s) ds \leq
2\mu^{4q_1+2}
K_1 \left(\frac{\mu n_t}{\l_t}\right)^{q_1}
\exp \left( -K_*\left(\frac{\mu n_t}{\l_t}\right)^{q_*} \right) \frac{\d(t)}{R_t}
\]
so
\be{Rdelta2}
\frac{R_t}{\d(t)}
\leq
2\mu^{4q_1+2}
K_1 \left(\frac{\mu n_t}{\l_t}\right)^{q_1}
\exp \left( -K_*\left(\frac{\mu n_t}{\l_t}\right)^{q_*} \right)
\ee
As we did in STEP 3, if $q_*,K_*$ are large enough, $R_{k}$ is small enough and the upper bound for the density holds on $\Delta_{k+1}$.
By using \eqref{fewR} and \eqref{Norm2}, we obtain
\[
\begin{split}
\Leb_n( y:|y-x_{k+1}|_{A_{R_{k+1}}(t_{k+1},x_{k+1})} \leq 1)
&\leq  2^n \Leb_n(y: |y-x_{k+1}|_{A_{R_k}(t_k,x_k)}\leq 1)  \\
&= 2^n
\sqrt{\det(A_{R_k} A_{R_k}^T(t_k,x_k))} \\
&=
C_k
\sqrt{\det(A A^T(t_k,x_k))}\,
R_k^{n-\frac{\dim\langle \s(t_k,x_k)\rangle}2},
\end{split}
\]
in which we have used the Cauchy-Binet formula (see also Remark \ref{C-B}).
Now, using the upper estimate for the density in the version of Theorem \ref{m} given in Remark \ref{est-tx}, we obtain
\[
\tilde{\PR}_k(X_{k+1}\in \Delta_{k+1})\leq
e^{\overline{C}_k}
\left(\frac{{R_k}}{\d_k}\right)^{n-\frac{\dim\langle \s(t_k,x_k)\rangle}2}
\]
where $\overline{C}_k=\overline{C}_{t_k}$, $\overline{C}\in \mathcal{A}$ (see the constant in the upper bound of \eqref{denco}).
Recall \eqref{Rdelta2}, for $t=t_k$
\[
\frac{R_k}{\d_k}
\leq
2\mu^{4q_1+2}
K_1 \left(\frac{\mu n_k}{\l_k}\right)^{q_1}
\exp \left( -K_*\left(\frac{\mu n_k}{\l_k}\right)^{q_*} \right)
\]
so we chose now $K_*,q_*$ large enough to have
\[
\tilde{\PR}_k(X_{k+1}\in\Delta_{k+1})\leq \exp(-K_2)
\]
for a constant $K_2>0$. From the definition of $N(T)$
\[
\int_0^T g_R(t)dt = \sum_{k=1}^{N(T)} \int_{t_{k-1}}^{t_k} g_R(t)dt = \frac{N(T)}{\mu^{2q_1+1}}\leq N(T).
\]
So, we have
\begin{align*}
& \PR\Big(\sup_{t\leq T} |X_t-x_t(\phi)|_{A_{R_t}(t,x_t(\phi))}\leq 1 \Big)\leq \E\Big(\prod_{k=1}^{N(T)} \tilde{\PR}_k(\Delta_{k+1})\Big)\\
& \quad \leq \prod_{k=1}^{N(T)} \exp(-K_2)= \exp(-K_2\,N(T))\leq
\exp\Big(-K_2\int_0^T g_R(t)\Big)
\end{align*}
and \eqref{tube-upper} holds.
\epr

\section{On the equivalence with the control distance}\label{sectioncontrol}

We establish here the local equivalence between the norm $%
\left\vert \cdot\right\vert _{A_{R}(t,x)}$ and the control (Carath\'e\-odory)
distance. We use in a crucial way the alternative characterization
given in \cite{NagelSteinWainger:85}. These results hold in the
homogeneous case, so we consider now the vector fields $\sigma _{j}(t,x)=\sigma _{j}(x)$, and the associated norm $A_{R}(t,x)=A_{R}(x)$. We assume in this section the following bound on $\s$:
there exists $\kappa\,:\,\R^n\to [1,+\infty)$ such that
\begin{equation}\label{bound}
\sup_{|y-x|\leq 1}\sum_{0\leq |\a| \leq 4}
\sum_{j=1}^d  |\partial_x^{\a} \s_j(y)|
\leq \kappa(x),
\quad\quad
\forall x \in \R^n.
\end{equation}
So, \eqref{bound} agrees with \eqref{Not2} in the homogeneous case and when $b=0$.

We now introduce a quasi-distance $d$ which is naturally associated to the family of norms $\left\vert y\right\vert _{A_{R}(x)}$.  We set
\[
O=\{x\in \R^n: \l_*(A(x))>0\}=\{x:\det(AA^T(x))\neq 0\}
\]
which is an open set since $x\mapsto \det(AA^T(x))$ is a continuous function. Notice that if $x\in O$ then $\det(A_R A_R^T(x))>0$ for every $R>0$.
For $x,y\in O$, we define $d(x,y)$ by
\[
d(x,y)<\sqrt{R}\quad \Leftrightarrow \quad \left\vert y-x\right\vert _{A_{R}(x)}<1.
\]
The motivation for taking $\sqrt{R}$ is the
following: in the elliptic case $\left\vert y-x\right\vert
_{A_{R}(x)}\sim R^{-1/2}\left\vert y-x\right\vert $ so $\left\vert
y-x\right\vert _{A_{R}(x)}\leq 1$ amounts to $\left\vert y-x\right\vert \leq
\sqrt{R}$.
It is straightforward to see that $d$ is a quasi-distance on $O$, meaning that $d$ verifies the following three properties (see \cite{NagelSteinWainger:85}):
\begin{itemize}
\item[$i)$] for every $x\in O$ and $r>0$, the set $\{y\in O\,:\,d(x,y)<r\}$ is open;

\item[$ii)$] $d(x,y)=0$ if and only if $x=y$;

\item[$iii)$] for every compact set $K\Subset O$ there exists $C>0$ such
that for every $x,y,z\in K$ one has $d(x,y)\leq C\big(d(x,z)+d(z,y)\big)$.
\end{itemize}
We recall the definition of equivalence of quasi-distances.
Two quasi-distances $d_{1}:\Omega \times\Omega \rightarrow \R^{+}$ and $%
d_{2}:\Omega\times \Omega\rightarrow \R^{+}$ are equivalent if for every
compact set $K\Subset\Omega$ there exists a constant $C$ such that for every
$x,y\in K$
\begin{equation}
\frac{1}{C}d_{1}(x,y)\leq d_{2}(x,y)\leq Cd_{1}(x,y).  \label{Norm7as}
\end{equation}
$d_1$ and $d_2$ are locally equivalent  if for every $\xi\in \Omega$ there exists a neighborhood $V$ of $\xi$ such that $d_1$ and $d_2$ are equivalent on $V$.

We introduce now the control metric. Without loss of generality, we assume $T=1$,

For $\psi\in L^2([0,1],\R^d)$, let $u(\psi)$ satisfy the following controlled equation:
\be{controlequation}
du_{t}(\psi )=\sum_{j=1}^{d}\sigma _{j}(u_{t}(\psi))\psi _{t}^{j}dt.
\ee
Notice that the equation for $u(\psi )$ is actually the skeleton equation \eqref{skeleton} when the drift $b$ is null.
For $x,y\in O$ we denote by $C^2_{\s,1}(x,y)$
the set of controls $\psi \in L^{2}([0,1];\R^d)$ such that the
corresponding solution $u(\psi)$ of \eqref{controlequation} satisfies $u_0(\psi)=x$ and  $u_{1}(\psi )=y$.   We define the control (Carath\'eodory) distance as
\begin{equation*}
d_{c}(x,y)=\inf_{\psi\in C^2_{\sigma,1}(x,y)}\|\psi\|_2.
% \Big\{\Big(\int_{0}^{1}\left\vert \psi _{s}\right\vert ^{2}ds%\Big)^{1/2}:\psi \in C_\s(x,y)\Big\}.
\end{equation*}

For $\d\in(0,1]$, we also denote $C^2_{\sigma,\delta }(x,y)$ the set of controls $\phi \in
L^{2}([0,\delta ];\R^d)$ such that the corresponding solution $u(\phi)$ to \eqref{controlequation} satisfies $u_{0}(\phi )=x$ and $u_{\delta }(\phi )=y$. For $\phi\in C_{\sigma,\delta}^2(x,y)$, we set the associated energy
\[
\varepsilon_{\phi }(\delta )=\Big(\int_{0}^{\delta }\left\vert
\phi _{s}\right\vert ^{2}ds\Big)^{1/2}.
\]
Notice that
\begin{equation}\label{dc-energy}
d_{c}(x,y)=
\sqrt{\delta } \inf_{\phi\in C^2_{\sigma,\delta}(x,y)}\varepsilon _{\phi }(\delta ).
\end{equation}
Indeed, for each $x,y\in\R^n$ and $\psi \in C^2_{\s,1}(x,y)$, take $\phi_t=\delta ^{-1}\psi (t\delta ^{-1})$ and $\xi_t=u_{t/\delta }(\psi )$. Then, $d\xi_{t}=\sum_{j=1}^{d}\sigma
_{j}(\xi_{t})\phi _{t}^{j}dt$ and of course $\xi_0=x$, $\xi_\delta=y$. Moreover, $\|\psi\|_2=\sqrt\delta\,\varepsilon_\phi(\delta)$.

Lastly, we define $C^\infty_{\sigma,1}(x,y)$ the set of paths $g \in L^\infty([0,1];\R^d)$ such that the corresponding solution $u(g)$ of \eqref{controlequation} satisfies $u_0(g)=x$ and $u_{1}(g)=y$. Using this set of controls, we define
\begin{equation*}
d_\infty (x,y)=\inf_{g\in C^\infty_{\sigma,1}(x,y)}\|g\|_\infty.
\end{equation*}

Under  \eqref{bound}, we define
\[
\D=\Big\{C:O\to \R_+\,:\,C=K \Big(\frac{\kappa(x)}{\lambda(x)}\Big)^{q},\ \exists\ K,q>0\Big\}.
\]
Notice that $\D$ is actually the set in \eqref{DD} in the homogeneous case.

\begin{theorem}

\label{NORM4}
Suppose that \eqref{bound} hold.

\smallskip

\textbf{A}. There exists $\bar{C}\in \D$ such that if $d_c(x,y)\leq 1/\bar{C}^2(x)$ then $d(x,y)\leq 2\bar{C}(x)d_c(x,y)$.

\smallskip

\textbf{B}. $d$ is locally equivalent to $d_{c}$ on $O$.

\smallskip

\textbf{C}. For every compact set $K\Subset O $ there exists $r_K$ and $C_K$ such that for every $x,y\leq r_K$ one has $d_c (x,y)\leq C_K d(x,y)$.

\end{theorem}

\bpr
\textbf{A.}
Assume that $d_c(x,y)\leq 1/\bar{C}^2(x)$, with $\bar C\in\mathcal{D}$ to be chosen later. We set $\delta(x)=\bar C^2(x) d_c(x,y)^2$. Notice that $\delta(x)\leq 1/\bar C^2(x)$. \eqref{dc-energy} with $\delta=\delta(x)$ gives
$$
d_{c}(x,y)=
\sqrt{\delta(x) } \inf_{\phi\in C^2_{\sigma,\delta(x)}(x,y)}\varepsilon _{\phi }(\delta(x))
=\bar{C}(x) d_c(x,y) \inf_{\phi\in C^2_{\sigma,\delta(x)}(x,y)}\varepsilon _{\phi }(\delta(x))
$$
and thus,
$$
\inf_{\phi\in C^2_{\sigma,\delta(x)}(x,y)}\varepsilon _{\phi }(\delta(x) )=\frac 1{\bar C(x)}<\frac 2{\bar C(x)}.
$$
Hence, there exists $\phi_*\in C^2_{\sigma,\delta(x)}(x,y)$ such that
$$
\varepsilon _{\phi_*}(\delta(x))<\frac 2{\bar C(x)}.
$$
For every fixed $x$, we apply Remark \ref{for-dc} to $\phi_*$ (recall that here $b\equiv 0$): there exists $\bar\delta,\bar\varepsilon\in 1/\mathcal{D}$ and $\bar C\in\mathcal{D}$ such that (with the slightly different notation of the present section)
$$
|u_{\delta}(\phi_*)-x|_{A_{\delta}(x)}\leq \overline{C}(x) (\ve_{\phi_*}(\delta)\vee \sqrt{\delta}),
$$
for every $\delta$ such that $\delta\leq \bar{\delta}(x)$ and $\varepsilon_{\phi_*}(\delta)\leq \bar{\varepsilon}(x)$. We have just proved that $\delta(x)\leq 1/\bar C^2(x)$ and $\varepsilon_{\phi_*}(\delta(x))\leq 2/\bar C(x)$. So, possibly taking $\bar C$ larger, we can actually use $\delta=\delta(x)$. And since $u_{\delta(x)}(\phi_*)=y$, the above inequality gives
$$
|y-x|_{A_{\delta}(x)}\leq \overline{C}(x) (\ve_{\phi_*}(\delta(x))\vee \sqrt{\delta(x)})\leq 2.
$$
By \eqref{Norm3}, we obtain $|y-x|_{A_{4\delta}(x)}\leq 1$, that is $d(x,y)\leq \sqrt{4\delta(x)}=2\bar C(x) d_c(x,y)$, and the statement follows.

\smallskip

\textbf{B.} We prove now the converse inequality. We use a result from
\cite{NagelSteinWainger:85}, for which we need to recall the definition of the quasi-distance $d_{\ast }$ (denoted by $\rho _{2}$ in \cite{NagelSteinWainger:85}). The definition we give here is slightly different but clearly equivalent.
For $\th \in \R^m$, consider the equation%
\be{Norm15}
dv_{t}(\th )=A(v_t(\th)) \th dt.
\ee
We denote
\[
\bar{C}_A(x,y)=\{\th\in \R^m\,:\, \mbox{ the solution $v(\theta)$ to \eqref{Norm15} satisfies $v_0(\theta)=x$ and $v_{1}(\th )=y$}\}.
\]
Notice that $\th\in\bar{C}_A(x,y)$ is a constant vector, and not a time depending control as in \eqref{controlequation}. Moreover, recalling the definitions \eqref{Al}-\eqref{Alucia} for $A$, \eqref{Norm15} involves also the vector fields $[\s_{i},\s_{j}]$, differently from  \eqref{controlequation}. In both equations the drift term $b$ does not appear.

Let $D_R$ be the diagonal matrix in \eqref{DR} and recall that $A_R(x)=A(x)D_R$.
We define
\[
d_{\ast }(x,y)=\inf\{R>0 \,:\, \mbox{ there exists $\th\in \bar{C}_A(x,y)$ such that $|D_R^{-1}\th|<1$}\}.
\]
As a consequence of Theorem 2 and Theorem 4 from \cite{NagelSteinWainger:85},
$d_{\ast }$ is locally equivalent with $d_\infty$. Since $d_c(x,y)\leq d_\infty(x,y)$ for every $x$ and $y$, one gets that $d_c$ is locally dominated from above by $d_{\ast}$.
To conclude we need to prove that $d_{\ast }$
is locally dominated from above by $d$.

Let us be more precise: for $x\in O$, we look for $C\in\D$ and $R\in 1/\D$ such that the following
holds: if $0<R \leq R(x)$ and $d(x,y)\leq \sqrt{R }$, then
there exists a control $\th\in \bar{C}_A(x,y)$ such that $|D_R^{-1}\th|<C(x)$. This implies $d_{\ast}(x,y)\leq C(x) \sqrt{R}$, and the statement holds. Notice that we discuss local equivalence, and that is why we can take $C(x)$ and $R(x)$ depending on $x$.

Recall that
$d(x,y)\leq \sqrt{R }$ means $|x-y|_{A_R(x)}\leq 1$, and this also implies $|x-y|\leq \l^*(A(x))\sqrt{R}$, by \eqref{Norm3}. Let $v(\theta)$ denote the solution to  \eqref{Norm15} with $v_0(\theta)=x$. We look for $\th$ such that  $v_1(\th)=y$. We define
\[
\Phi(\th)= \int_{0}^{1} A(v_s(\th)) \th ds = A(x) \th + r(\th)
\]
with $r(\th)=\int_{0}^{1} (A(v_s(\th))-A(x))\th ds$. With this notation, we look for $\th$ such that $\Phi(\th)=y-x$. We introduce now the Moore-Penrose pseudoinverse of $A(x)$: $A(x)^+=A(x)^T(A A^T(x))^{-1}$. The idea here is to use it as in the least squares problem, but we need some computations to overcome the fact that we are in a non-linear setting. We use the following properties: $A A(x)^+ = \Id$; $|x-y|_{A(x)}=|A(x)^+\,(x-y)|$. Write $\th=A(x)^+ \gamma$, $\gamma\in \R^n$. This implies $A(x)\th =\gamma$, and so we are looking for $\gamma\in\R^n$ such that
\[
\gamma + r( A(x)^+ \gamma)=y-x.
\]
One has $r(0)=0$, $\nabla r(0)=0$ and, as a consequence, $|r(\th)|\leq C(x) |\th|^2$, for some $C\in\mathcal{D}$ -- from now on, $C\in \D$ will denote a function that may vary from line to line.

From the local inversion theorem (in a quantitative form), there exists $l\in \D$ such that $\gamma\mapsto\gamma + r( A(x)^+ \gamma)$ is a diffeomorrphism from $B(0,l_x)$ to $B(0,l_x/2)$. Remark that $|x-y|\leq \l^*(A(x))\sqrt{R}$. So, taking $R_x$ such that $\l^*(A(x))\sqrt{R}=l_x/2$, then for every $R<R_x$ and $|y-x|<\lambda^*(A(x))\sqrt R$ then there exists a unique $\gamma$ such that $\gamma + r( A(x)^+ \gamma)=y-x$ and moreover, $|\gamma|\leq 2|x-y|$.
Now, using \eqref{Norm3}
\[
|r( A(x)^+ \gamma)|_{A_R(x)}\leq \frac{\l^*(A(x)) |r( A(x)^+ \gamma)|}{R}
\leq C_x\frac{|A(x)^+ \gamma|^2}{R}\leq C_x\frac{|x-y|^2}{R} \leq C_x |x-y|_{A_R(x)}^2.
\]
Since $\gamma =x-y-r(A(x)^+ \gamma)$,
\[
|\gamma|_{A_R(x)}\leq |x-y|_{A_R(x)}+C_x|x-y|_{A_R(x)}^2\leq C_x,
\]
(using $|x-y|_{A_R(x)}\leq 1$). We have
$|D_R^{-1}\th|=|D_R^{-1} A(x)^+ \g|$.
Since $A_R^+ A_R(x)=A_R^+(x) A(x) D_R$ is an orthogonal projection and $AA^+(x)$ is the identity,
\[
|D_R^{-1}\th|\leq |D_R^{-1} A(x)^+ \g|\leq A
|A_R^+(x) A(x) D_R \,D_R^{-1} A(x)^+ \g|
=|A_R^+(x) \g|=
|\g|_{A_R(x)}.
\]
So $|D_R^{-1}\th|\leq C_x$, and this implies $d_*(x,y)\leq C_x \sqrt{R}$.

\smallskip

\textbf{C.}
The proof immediately follows from the previous items.
\epr

The proof of Theorem \ref{tubecontrold} is now an immediate consequence of Theorem \ref{mttubesS} and Theorem \ref{NORM4}. The only apparent problem is that in Theorem
\ref{NORM4} the global estimate \eqref{bound} is required, whereas in Theorem \ref{tubecontrold} the local estimate $(H_1)$ in \eqref{Not6} holds. But this is not really a problem, since it can be handled as already done for Theorem \ref{tubecontrold} (see Remark \ref{2}).

\bibliographystyle{plain}
\bibliography{bibliografia}

\end{document}